\documentclass[11pt]{article}
\usepackage{graphicx}
\usepackage{color}
\usepackage{amsmath,amsfonts,amssymb,amsthm}
\def\ds{\displaystyle}
\def\forall{\hbox{for all}~}
\def\L{{\bf L}}

\def\ve{\varepsilon}

\def\R{{\mathbb R}}

\def\vp{\varphi}

\def\tv{\hbox{Tot.Var.}}
\def\vs{\vskip 2em}

\def\v{\vskip 1em}
\def\D{{\cal D}}
\def\O{{\cal O}}
\def\N{{\cal N}}

\def\begi{\begin{itemize}}
\def\endi{\end{itemize}}
\def\C{{\cal C}}
\def\dint{\int\!\!\!\int}

\def\Tilde{\widetilde}

\def\bega{\begin{array}}
\def\enda{\end{array}}
\def\meas{\hbox{meas}}
\def\bel{\begin{equation}\label}
\def\eeq{\end{equation}}
\def\sqr#1#2{\vbox{\hrule height .#2pt
\hbox{\vrule width .#2pt height #1pt \kern #1pt
\vrule width .#2pt}\hrule height .#2pt }}
\def\square{\sqr74}
\def\endproof{\hphantom{MM}\hfill\llap{$\square$}\goodbreak}

\newtheorem{theorem}{Theorem}[section]
\newtheorem{corollary}{Corollary}[section]
\newtheorem{lemma}{Lemma}[section]

\newtheorem{remark}{Remark}[section]
\newtheorem{definition}{Definition}[section]
\newtheorem{example}{Example}[section]

\begin{document}
\title{\bf  %Uniqueness and Approximation Estimates for 
One Dimensional Hyperbolic Conservation Laws: Past and Future}
\vs
%\author{Alberto Bressan~(review paper)}
\author{Alberto Bressan\\
~~~\\
Department of Mathematics, Penn State University, \\
University Park, Pa.~16802, USA.\\
bressan@math.psu.edu~}
\maketitle

\begin{abstract} Aim of these notes is to provide a brief review of the current well-posedness 
theory for hyperbolic systems of conservation laws in one space dimension, also pointing out 
open problems and possible research directions.
They supplement
the slides of the short course given by the author in Erice, May 2023, available at:
sites.google.com/view/erice23/speakers-and-slides. \end{abstract}

\section{Introduction}
\label{sec:1}
\setcounter{equation}{0}
Aim
of these notes is provide a brief review of the current well-posedness theory for hyperbolic systems of conservation laws in one space dimension, also pointing out 
open problems and possible research directions.
They supplement
the slides of the short course given by the author in Erice, May 2023, available at:
sites.google.com/view/erice23/speakers-and-slides.

Section~2 introduces basic definitions, 
including the concept of weak solution, and various admissibility conditions.
Section~3  describes several approximation methods.
The main results on global existence of weak solutions and their continuous dependence on initial data are recalled in Sections 4 and 5.  The recent advances, on the
uniqueness of weak solutions that satisfy the Liu admissibility condition,
are covered in greater detail in Section~6.  The relevance of these results 
toward error
bounds for all kinds of approximate solutions is discussed in Section~7, 
together with two specific open problems.  
Finally, Section~8 is devoted to solutions with 
possibly unbounded total variation, recalling the main known results and pointing out
some research directions.

 \v
\section{Basic concepts}
\label{sec:2}
\setcounter{equation}{0}

%%%%%
\subsection{Hyperbolic systems.}
A system of conservation laws in one space dimension has the form
\bel{hcl} u_t+f(u)_x~=~0.\eeq
In components, this can be written as 
\bel{2}
{\partial\over\partial t}\begin{pmatrix} u_1\cr \vdots \cr u_n\end{pmatrix}
+{\partial\over\partial x}\begin{pmatrix} f_1(u)\cr \vdots \cr f_n(u)\end{pmatrix}
~=~\begin{pmatrix} 0\cr \vdots \cr 0\end{pmatrix}.\eeq
Here $u = (u_1,\ldots, u_n)^T$ is the vector of {\bf conserved quantities}, 
while $f=(f_1,\ldots, f_n)^T$ is the vector of  {\bf fluxes}.
Conservation laws provide the fundamental mathematical models in continuum physics~\cite{Dbook}.
A primary example is provided by the Euler equations of gas dynamics, 
accounting for the conservation of mass, momentum and energy~\cite{Euler}.

%
%
%\begin{figure}[htbp]
%\centering
%  \includegraphics[scale=0.5]{FIG/claw16.eps}
%    \caption{\small  Deriving the system of equations (\ref{hcl}) from the conservation 
%    relations.}
%\label{f:claw16}
%\end{figure}
%

Smooth solutions to the system of PDEs (\ref{hcl}) can be obtained by
solving the equivalent
quasilinear system
\bel{qls}
u_t + A(u) u_x~=~0,\qquad \qquad \hbox{where}\quad A(u) \doteq  Df(u).\eeq
We say that the system (\ref{hcl}) is {\bf strictly hyperbolic} if at every point $u$
the $n\times n$ Jacobian matrix $A(u)\doteq Df(u)$ has $n$ real distinct eigenvalues
\bel{eval} \lambda_1(u)~<~\lambda_2(u)~<~\cdots~<~\lambda_n(u).\eeq
When this holds, one can find bases of right and left eigenvectors, say
$\{r_1,\ldots, r_n\}$, $\{l_1,\ldots, l_n\}$, with
\bel{evect}A(u) r_i(u)\,=\,\lambda_i(u) r_i(u),\qquad\qquad l_i(u) A(u) \,=\,\lambda_i(u)l_i(u),
\qquad\qquad i=1,\ldots,n.\eeq
These vectors can be normalized so that
$$\bigl| r_i(u)\bigr|\,=\,1,\qquad\qquad l_i(u) r_j(u)~=~\left\{\bega{rl} 1\quad &\hbox{if}~~i=j,\cr
 0\quad &\hbox{if}~~i\not= j.\enda\right.$$
The behavior of eigenvalues of $Df(u)$ strongly affect the nature of solutions to
(\ref{hcl}).    Following classical literature \cite{Lax}, we say that the $i$-th characteristic
field is {\bf genuinely nonlinear} if the directional derivative of the eigenvalue 
$\lambda_i$ in the direction of the corresponding eigenvector $r_i(u)$ satisfies
\bel{gnl} \nabla \lambda_i(u)\cdot r_i(u)~>~0\qquad\qquad\forall u.\eeq
On the other hand, we say that the $i$-th characteristic
field is {\bf linearly degenerate} if 
\bel{ldeg} \nabla \lambda_i(u)\cdot r_i(u)~=~0\qquad\qquad\forall u.\eeq
Throughout the following we assume that the flux function $f$ is 
at least twice continuously 
differentiable, so that the above derivatives are well defined.

\begin{example}\label{ex:1} {\rm 
In Lagrangian coordinates, the Euler equations of {\bf isentropic gas dynamics}
take the form
\bel{igd} \left\{ \begin{array}{cl} v_t-u_x&=~0,\cr
u_t+p(v)_x&=~0.\end{array}\right.\eeq
 Here $\rho$ is the density of the gas, $v=\rho^{-1}$ is specific volume,  $u$ is the  velocity
and  $p=p(v)$ is the pressure.   A natural choice for the pressure is
$p(v)=k v^{-\gamma}$,
with $1\leq \gamma\leq 3$.
The eigenvalues of the Jacobian matrix 
$$A\,\doteq\, Df\,=\,\begin{pmatrix}0 & -1 \cr p'(v) & 0\cr\end{pmatrix}$$ are
$$\lambda_1~=~ -\sqrt{-p'(v)}\,,\qquad\qquad
\lambda_2~= ~\sqrt{-p'(v)}$$
Since
$p'(v)<0$, the system is strictly hyperbolic.   A further computation reveals that 
both characteristic fields are genuinely nonlinear.}
\end{example}
%\v
%The {\bf Cauchy problem} for the system (\ref{hcl}) consists in finding a solution
%$u=u(t,x)$ defined for all $t\geq 0$, $x\in\R$, with the initial data
%\bel{id}u(0,x)~=~\bar u(x),\qquad\qquad x\in\R.\eeq

\subsection{Weak solutions.}
\begin{figure}[htbp]
\centering
  \includegraphics[scale=0.5]{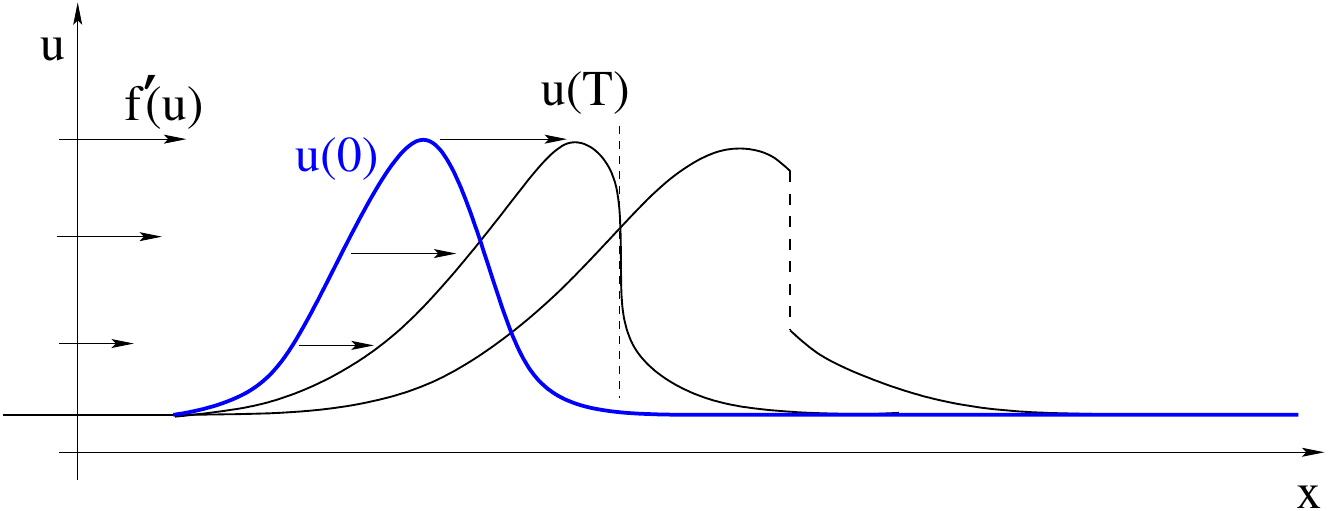}
    \caption{\small An example where the gradient 
    $u_x$ of the solution  becomes unbounded at a finite time $T$. For $t>T$, the
    solution contains a shock and must be interpreted in distributional sense.  }
\label{f:claw18}
\end{figure}

A key feature  of hyperbolic conservation laws is that, even for initial data
$\bar u\in \C^1$ (i.e., continuously differentiable), the gradient $u_x$ of the
solution may blow up at a finite time $T$  (see Fig.~\ref{f:claw18}).
In order to prolong the solution  also for $t>T$, one must work within a space
of discontinuous functions, interpreting the  equation (\ref{hcl}) in distributional sense.
In the following, $\L^1_{loc}(\Omega)$ denotes the space of locally integrable functions
defined on an open subset $\Omega\subset\R\times\R$, with values in $\R^n$.   
Moreover, 
$\C^1_c(\Omega)$ denotes the space of continuously differentiable 
functions with compact support. 

\begin{definition}\label{d:ws1} 
Let 
$u = u(t,x)$ be a function defined on an open set $\Omega\subseteq\R\times\R$.
We say that $u$ is a {\bf weak solution} to the system of conservation laws (\ref{hcl})
if $u, f(u)\in \L^1_{loc}(\Omega) $ and 
\bel{weaksol}\dint\big\{
u\phi_t+f(u)\phi_x\big\}~dxdt~=~0\qquad
\hbox{for all} ~~~\phi\in \C^1_c(\Omega).\eeq
\end{definition}

\v
\begin{definition}\label{d:ws2}  A function $u = u(t,x)$ is a {\bf  weak solution}  to the {\bf Cauchy problem}
\bel{CP}u_t + f(u)_x~=~0,\qquad\qquad u(0,x)=\bar u(x),\qquad\qquad t\in [0,T],\eeq
if the map $t\mapsto u(t,\cdot)$ is continuous with values in $\L^1(\R;\,\R^n)$, satisfies the initial 
condition in (\ref{CP}), and moreover
$$\int_0^T\!\int_{-\infty}^{+\infty} \big\{
u\phi_t+f(u)\phi_x\big\}~dxdt ~=~0\qquad
\hbox{for all} ~~~\phi\in \C^1_c\bigl(]0,T[\times\R\bigr)$$
\end{definition}

\begin{figure}[htbp]
\centering
  \includegraphics[scale=0.5]{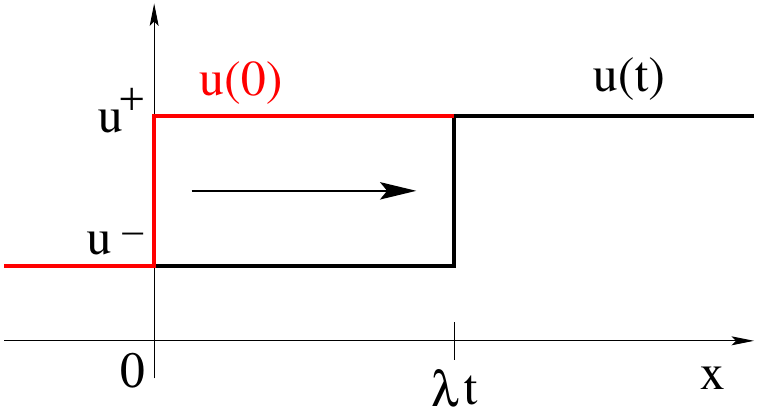}
    \caption{\small A shock with left and right states $u^-, u^+$, moving with speed $\lambda$.  }
\label{f:claw20}
\end{figure}

Notice that the identity (\ref{weaksol}) is obtained  multiplying (\ref{hcl}) by the 
test function $\phi$ and integrating by parts.

The simplest example of a discontinuous solution to (\ref{hcl}) is a single shock,
shown in Fig.~\ref{f:claw20}.
\bel{shock} u(t,x)~=~\left\{\bega{rl}u^-\quad &\hbox{if}\quad x<\lambda t\,,\cr
u^+\quad &\hbox{if}\quad x>\lambda t\,.\enda\right.\eeq
It is well known that the above function is a weak solution 
if and only if the shock speed $\lambda$ and the left and right states $u^-, u^+$ 
satisfy the {\bf Rankine-Hugoniot equations}
\bel{RH1}\lambda\cdot (u^+-u^-)~=~f(u^+)-f(u^-).\eeq
In other words, the vector equation (\ref{RH1}) states that
$$\hbox{[speed] $\times$ [jump in the state] ~=~ [jump in the flux].}$$

\begin{figure}[htbp]
\centering
  \includegraphics[scale=0.5]{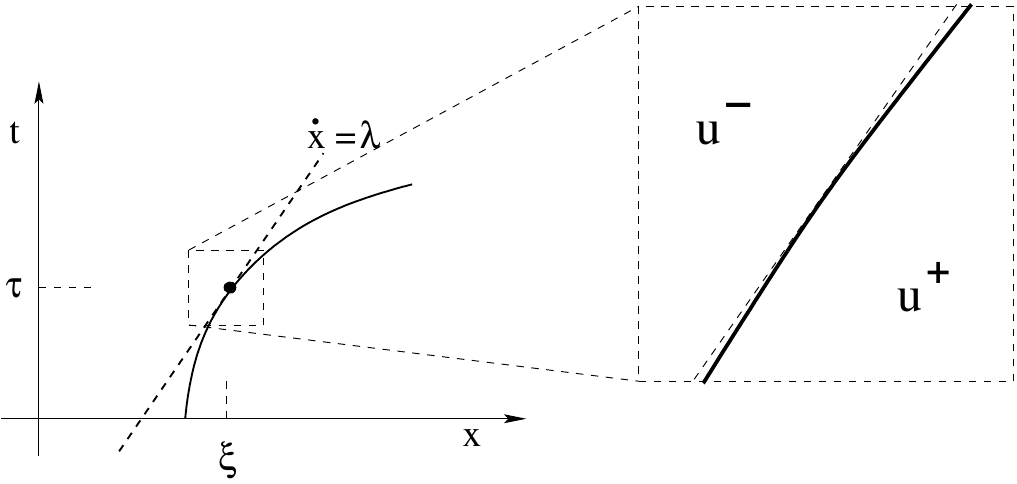}
    \caption{\small  A point of approximate jump.    }
\label{f:f146}
\end{figure}
To state a version of the Rankine-Hugoniot conditions which applies to more general solutions,
we introduce
\begin{definition}\label{d:ajump}
The function $u=u(t,x)$ has an {\bf approximate jump}
at the point $(\tau,\xi)\in \R^2$ if there exists vectors $u^+\not= u^-$ and a speed $\lambda$ such that, setting
\bel{Udef}U(t,x)~\doteq~\left\{\bega{rl}  u^-\qquad &\hbox{if}\qquad x<\lambda\,t,\cr
 u^+\qquad &\hbox{if}\qquad x>\lambda\,t,\enda\right.
\eeq
one has
\bel{ajump} \lim_{r\to 0+} \,{1\over r^2}
\int_{-r}^r\int_{-r}^r \Big|
u(\tau+t,\,\xi+x)-U(t,x)\Big|\,dxdt~=~0.\eeq
We
say that $u$ is {\bf approximately continuous} at the point
$(\tau,\xi)$ if (\ref{ajump}) holds with $u^+=u^-$ (and
$\lambda$ arbitrary).
\end{definition}

\begin{theorem} {\bf (Rankine-Hugoniot equations)}  Let
$u$ be a bounded weak solution of (\ref{hcl}),
having an
approximate jump at a point $(\tau,\xi)$.
Then the left and right states $u^-, u^+$ and the speed $\lambda$ satisfy the Rankine-Hugoniot equations (\ref{RH1}).\end{theorem}

For a proof, see for example~\cite{Bbook}.
Writing the Rankine-Hugoniot equations in the form
$$\bega{rl}
\lambda\,
(u^+-u^-)&=~f(u^+)-f(u^-)~=~\displaystyle\int_0^1 Df\bigl( \theta
u^++(1-\theta)u^-\bigr)\cdot (u^+-u^-)~d\theta\\[4mm]
 &=~  A(u^+,u^-)\cdot
(u^+-u^-),\end{array}$$
we see that
\begi
\item[(i)]  The jump $u^+-u^-$ is an eigenvector of the
averaged matrix $A(u^+,u^-)$.

\item[(ii)] The speed $\lambda$ coincides
with the corresponding eigenvalue.
\endi

\subsection{Admissibility conditions.}
In general, solutions to the Cauchy problem (\ref{CP}) may not be unique,
as soon as discontinuities are present.
To single out a unique weak solution one needs to impose further admissibility conditions
on the shocks.  These can be derived by three different approaches:
\begi
\item 
Stability w.r.t.~small perturbations.
\item Vanishing viscosity approximations.
\item Entropy dissipation.
\endi

\begin{figure}[htbp]
\centering
  \includegraphics[scale=0.45]{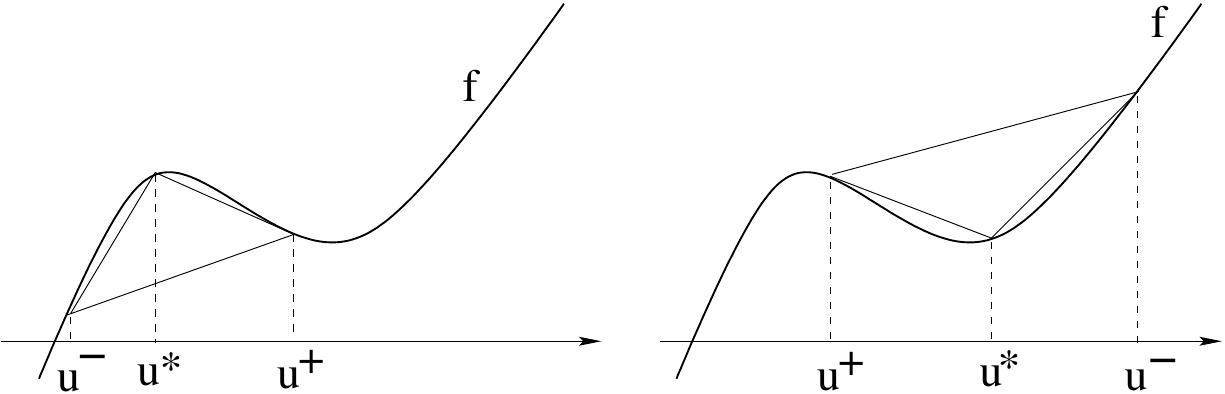}
    \caption{\small For a scalar conservation law,  according to  (\ref{RH2}) 
     the speed of a shock
    is the slope of a secant line to the graph of $f$. }
\label{f:f54}
\end{figure}

{\bf 1. A stability condition.}
Consider first a scalar conservation law. In this case, the Rankine-Hugoniot condition (\ref{RH1})
simply states that the speed of a shock with left and right states $u^-, u^+$ must be
\bel{RH2}\lambda~=~{f(u^+)-f(u^-)\over u^+-u^-}\,.\eeq
Looking at the graph of the function $f(u)$ (see Fig.~\ref{f:f54}), this means
$$\hbox{speed of the shock~~=~~slope of the secant line through $u^-, u^+$.}$$

To check the stability of the solution (\ref{shock}),
we can perturb the shock 
by inserting an intermediate state $u^*\in [u^-, u^+]$

\begin{figure}[htbp]
\centering
  \includegraphics[scale=0.4]{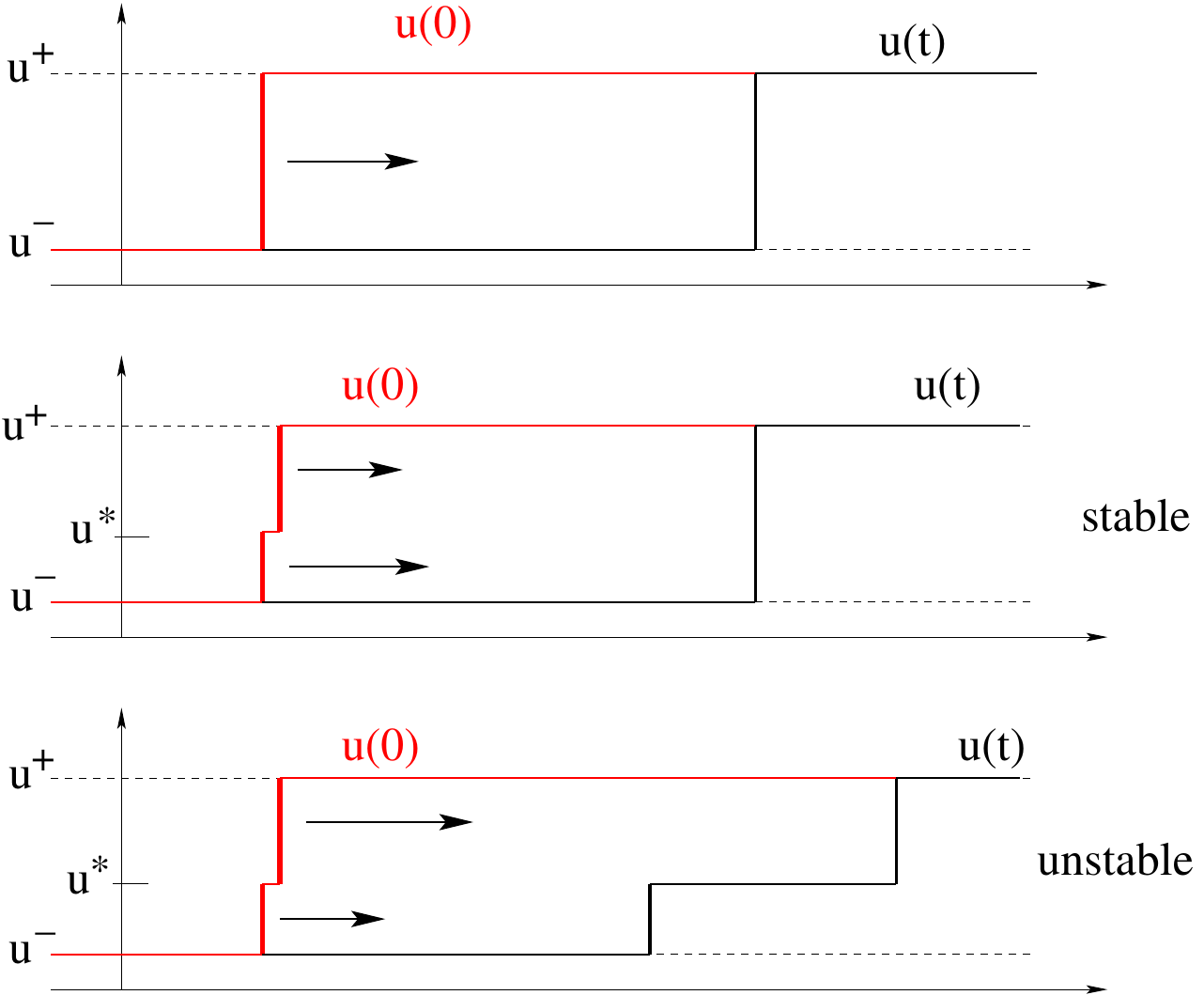}
    \caption{\small  A solution containing a single shock (top figure) can be perturbed 
    into a solution which initially contains two nearby shocks.   If the shock behind travels faster than the shock ahead (center figure), then the original shock is stable.   If the shock behind is slower  than the shock ahead (lower figure),  the original shock is unstable.}
\label{f:claw22}
\end{figure}

The original  shock will be stable w.r.t.~this perturbation iff
$$\hbox{[speed of jump behind] ~$\geq$ ~[speed of jump ahead]}.$$
By (\ref{RH2}), this means
\bel{stab1}{f(u^*)-f(u^-)\over u^*-u^-}~\geq~
{f(u^+)-f(u^*)\over u^+-u^*}\,.\eeq
Interpreting the two sides of (\ref{stab1}) as slopes of secant lines to the graph of $f$,
as shown in Fig.~\ref{f:f54} one obtains the following {\bf stability conditions.}
\begi
\item[(i)] when $u^-<u^+$ the graph of $f$ should remain above
the secant line through $u^-, u^+$.
\item[(ii)] when $u^->u^+$, the graph of $f$ should
remain below the secant line through $u^-, u^+$.
\endi
The two above cases are equivalent to one single inequality, namely
\bel{ac4} \hbox{speed of the shock}~ [u^-, u^+]~~\leq~~\hbox{speed of any intermediate shock}~ [u^-, u^*].\eeq
For every intermediate state $u^*$ between $u^-$ and $u^+$, the stability of the shock thus requires
\bel{stab2}
{f(u^+)-f(u^-)\over u^+-u^-}~\leq~{f(u^*)-f(u^-)\over u^*-u^-}\,.\eeq

\begin{figure}[htbp]
\centering
  \includegraphics[scale=0.45]{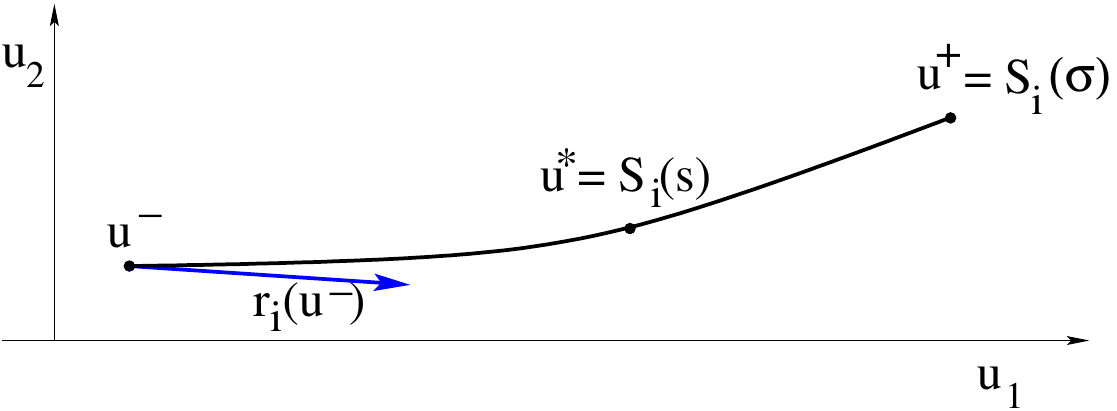}
    \caption{\small The $i$-shock curve $s\mapsto S_i(s)$ through the left state $u^-$.}
\label{f:hyp253}
\end{figure}

The formulation (\ref{stab2}) is particularly important, because it can be extended to any $n\times n$ strictly hyperbolic system of conservation laws.  
To state this general admissibility condition,
we recall that, for any given left state $u^-\in \R^n$ and $i\in \{1,\ldots,n\}$, 
one can find a curve
$s\mapsto S_i(s)$ of right states which can be connected to $u^-$ by an $i$-shock.
More precisely (see Fig.~\ref{f:hyp253}),
$$S_i(0) \,=\, u^-,\qquad\qquad {d\over ds} S_i(s)\bigg|_{s=0} ~=~r_i(u^-),$$
and for every $s$ the Rankine-Hugoniot equations hold:
\bel{RH3} f\bigl( S_i(s)\bigr) - f(u^-)~=~\lambda_i(s) \bigl( S_i(s)-u^-\bigr),\eeq
for some speed $\lambda_i(s)$.

The Liu admissibility condition for general shocks can now be stated as follows.
\begin{definition}\label{d:liu}
A shock of the $i$-th family, connecting the states $u^-$ and $u^+= S_i(\sigma)$
is {\bf  Liu-admissible}  if the speeds of all intermediate shocks satisfy
\bel{liuadm}\lambda_i(s)~\geq~\lambda_i(\sigma)\qquad\qquad
\hbox{for all}~~~s\in[0,\sigma].\eeq
\end{definition}

\begin{figure}[htbp]
\centering
  \includegraphics[scale=0.35]{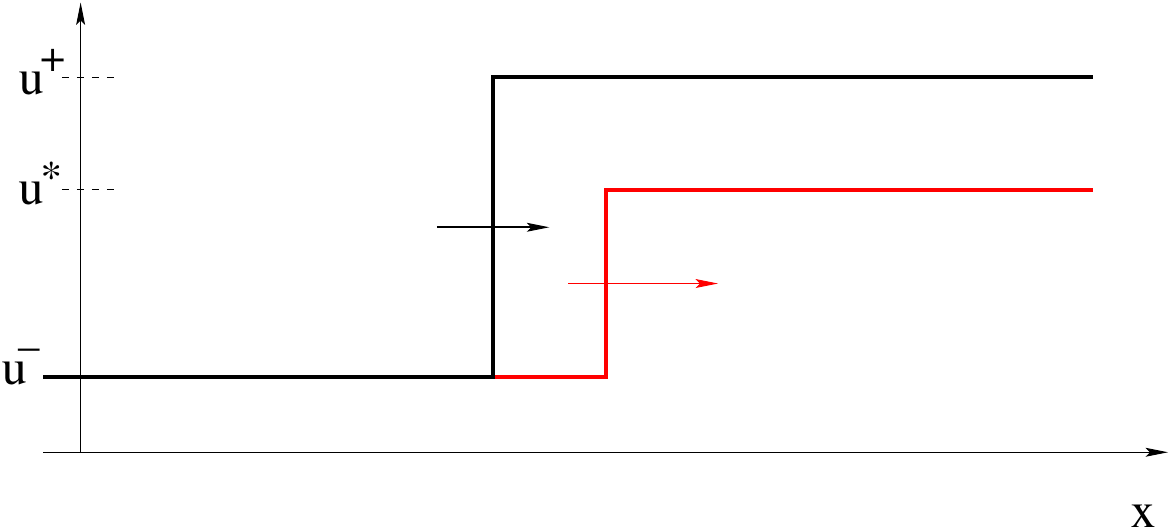}
    \caption{\small According to the Liu condition, the shock with left and right states
    $u^-, u^+$ is admissible
    if its speed is slower that the speed of every intermediate shock, joining the states
    $u^-, u^*$. }
\label{f:hyp257}
\end{figure}

{\bf 2. Vanishing viscosity limits.}  From physical considerations, 
it is often natural to assume that
the  ``good"  solutions to the system of conservation laws (\ref{hcl}) 
are those obtained as limit of solutions to vanishing viscosity approximations
\bel{vva}u^\ve_t + f(u^\ve)_x~=~\ve\,u^\ve_{xx}\,.\eeq
In particular, we say that the shock solution (\ref{shock}) is 
{\bf vanishing viscosity admissible} if
it can be obtained as a limit of solutions to (\ref{vva}). 

In the case of a single shock,
the equivalence between the  vanishing viscosity and the  Liu admissibility condition 
for a single shock has been proved in  \cite{Liu76} and in \cite{Bi}.  
More generally, the analysis in 
\cite{BiB} shows that every limit of vanishing viscosity approximations satisfies the 
Liu condition at every point of approximate jump.
\v
{\bf 3. Entropy admissibility condition.} 
Given the hyperbolic system of conservation laws (\ref{hcl}), 
a scalar function $\eta(u)$ is called an {\bf entropy}  with {\bf  entropy flux} $q(u)$ if
\bel{eef} D\eta(u)\cdot Df(u)~=~Dq(u).\eeq
We observe that (\ref{eef}) is a system of
$n$ first order PDEs  for the 
2  functions $\eta,q$ of the variables $(u_1,\ldots, u_n)$.   In general, this is
overdetermined and has no solution if $n>2$.   
However, there are relevant physical systems
where a nontrivial entropy can still be found.

By (\ref{eef}),
every smooth solution to (\ref{hcl}) satisfies the additional conservation law
$$\eta(u)_t + q(u)_x~=~D\eta(u)\,u_t +D\eta(u) Df(u)u_x~=~0.$$
On the other hand, the entropy may not be conserved in the presence of shocks.
\v

\begin{definition} Assume that the 
 hyperbolic system of conservation laws (\ref{hcl}) admits a {\bf  convex entropy} $\eta(u)$  with {\bf 
entropy flux} $q(u)$.
We say that a weak solution $u=u(t,x)$ is {\bf entropy admissible} if it satisfies
the inequality 
\bel{eac}\eta(u)_t + q(u)_x~\leq~0\eeq
in distributional sense. That means
\bel{eac2} \dint\bigl\{\eta(u) \phi_t+ q(u)\phi_x\bigr\}\, dxdt~\geq~0
\qquad\forall \phi\in \C^1_c\,,~~~\vp\geq 0.\eeq
\end{definition}
As a special case, the shock solution  at  (\ref{shock}) is entropy-admissible iff
\bel{shad}
\lambda\bigl[ \eta(u^+)-\eta(u^-)\bigr]~\geq~ q(u^+)-q(u^-).\eeq

It is well known that, if a convex entropy exists, 
every limit of vanishing viscosity approximations satisfies the 
entropy admissibility conditions \cite{Bbook, Dbook}.

\begin{remark}{\rm  In the classical theory of gas dynamics, the second law of thermodynamics implies that the physical entropy should increase in time.
To reconcile this fact with the decrease in the entropy stated at (\ref{eac}), it suffices to observe that the physical entropy is concave down, while the entropies considered by the mathematical theory are always 
 convex functions. Hence the change in the sign.}
\end{remark}

\section{Approximation methods}
\label{sec:3}
\setcounter{equation}{0}
Several techniques for constructing approximate solutions
to the Cauchy problem
\bel{CP3} u_t + f(u)_x~=~0,\qquad \qquad u(0,\cdot) = \bar u,\eeq
have been considered in the literature.  Generally speaking, these methods are known
to converge to the exact solution in two main cases:
\begi
\item[(i)] For a scalar conservation law, based on Kuznetsov's estimates. See for example
\cite{HR}.
\item[(ii)] For general $n\times n$ systems, as long as the exact solution remains $\C^1$, i.e., continuously differentiable.
\endi

On the other hand, for weak solutions to $n\times n$ hyperbolic systems,
the convergence of approximations requires a careful analysis. 
In various cases, the convergence still remains an open problem.

\begin{figure}[htbp]
\centering
  \includegraphics[scale=0.45]{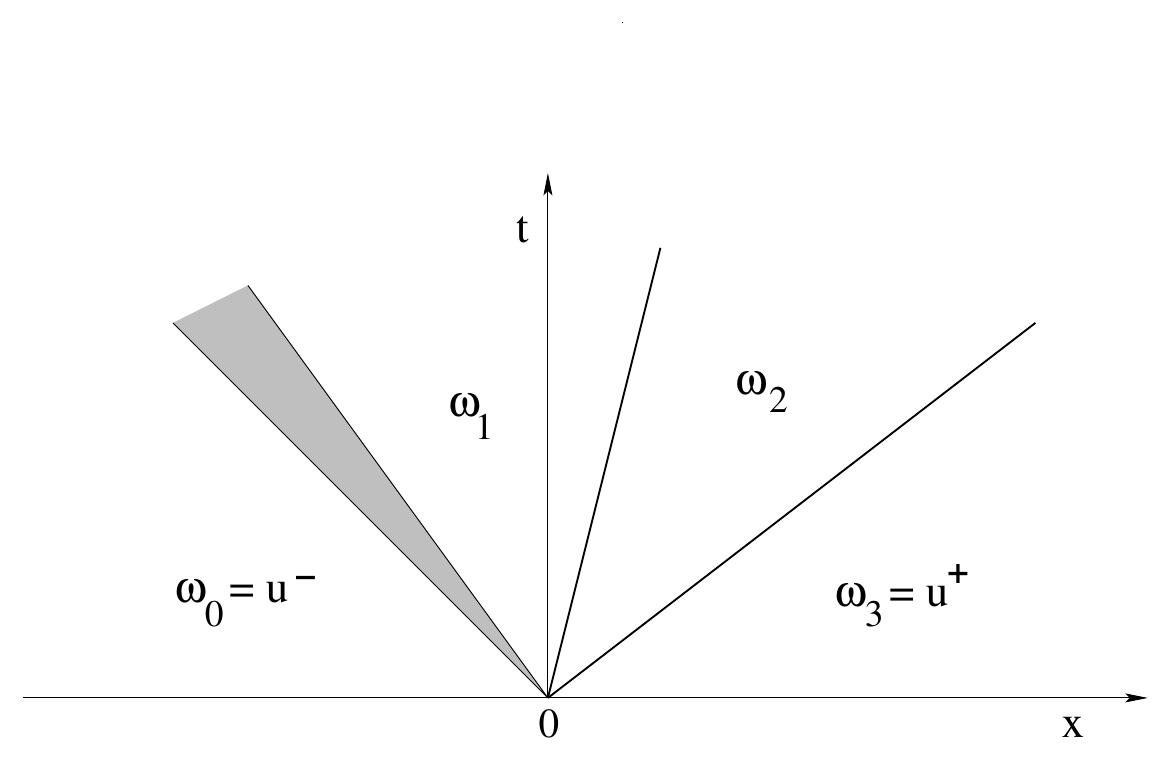}
    \caption{\small A typical solution to a Riemann problem for a $3\times 3$ hyperbolic system.}
\label{f:f136}
\end{figure}

In many algorithms,
a basic building block is provided by the {\bf Riemann problem}, where the  initial datum
is piecewise constant  with one single jump at the origin:
\bel{RP}
u(0,x)~=~\left\{ \bega{cl} u^-\quad &\hbox{if}~~x<0,\cr
u^+\quad &\hbox{if}~~x>0.\enda\right.\eeq
For an $n\times n$ hyperbolic system, assuming that every characteristic field
is either genuinely nonlinear or linearly degenerate, the general solution to the Riemann 
problem was 
first constructed by Lax~\cite{Lax}.  As shown in Fig.~\ref{f:f136},
it consists of $n+1$ constant states
$$u^-=\omega_0\,, ~~\omega_1\,,\qquad \cdots\qquad ,~\omega_n=u^+,$$
where each couple of states $(\omega_{k-1}, \, \omega_k)$ are separated by an 
admissible shock
or by a centered rarefaction wave of the $k$-th family.
Solutions $u=u(t,x)$  to a Riemann problem are self-similar, in the sense that
they are invariant w.r.t.~a rescaling symmetry:
$$u(t,x)~=~u(\theta t,\,\theta x)\qquad \forall \theta>0.$$

A brief survey of different approximation methods for the Cauchy problem (\ref{CP3}) is given below.
\v
{\bf 1. The upwind  Godunov scheme.} 
To simplify the presentation, we shall assume that 
all characteristic speeds satisfy 
\bel{la1}\lambda_i(u)\in [0,1].\eeq
This is not restrictive, because if $\lambda_i(u)\in [-M,\, M]$ one can  achieve (\ref{la1})
by the simple coordinate change
\bel{cchange}\tilde x\,=\,x+Mt,\qquad\qquad \tilde t\,=\,2Mt.\eeq

\begin{figure}[htbp]
\centering
  \includegraphics[scale=0.4]{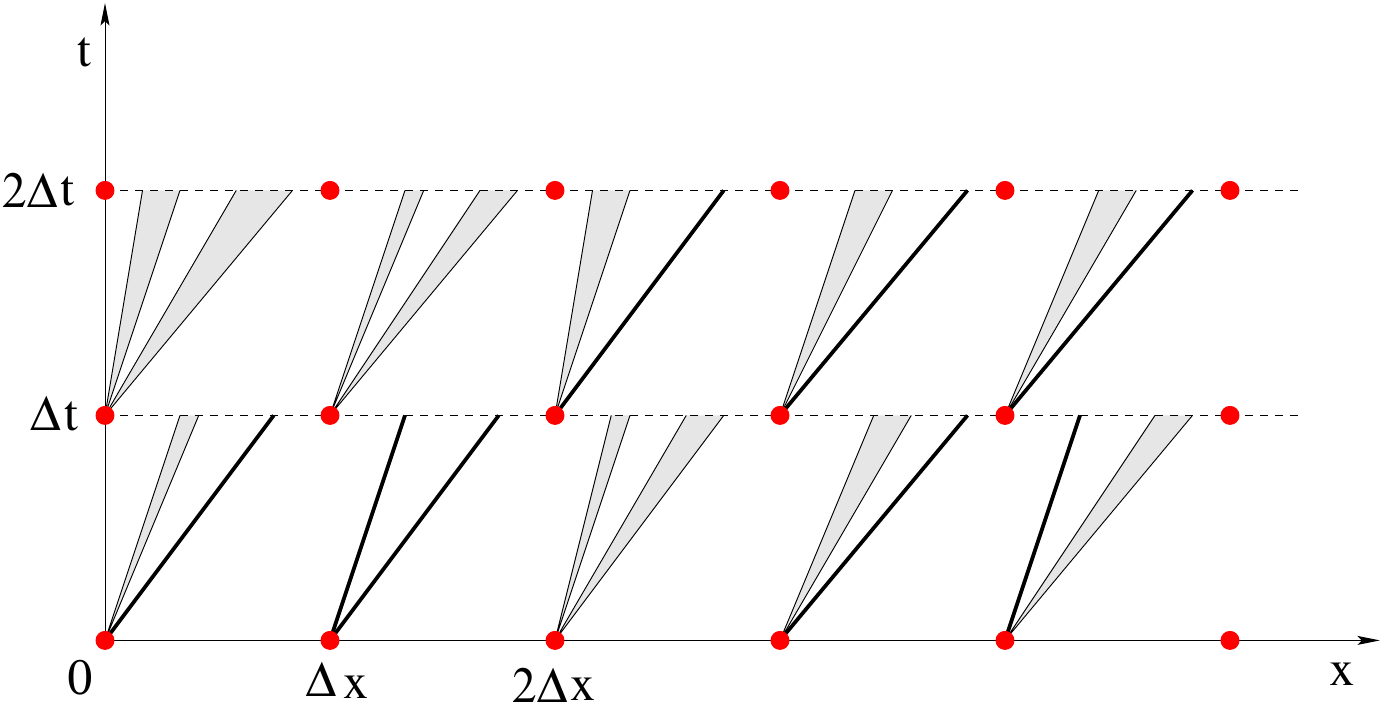}
    \caption{\small An approximate solution obtained by solving a  Riemann problem at each node of the grid.
    }
\label{f:claw46}
\end{figure}

The Godunov (upwind) scheme starts by constructing a grid in the 
$t$-$x$ plane with step size $\Delta t = \Delta x=\ve$, see Fig.~\ref{f:claw46}.
\begi
\item The grid points are  $(t_j, x_k) = (j\cdot \Delta t\,,~ k\cdot \Delta x)$.
%\item At time $t_0=0$ the initial datum $\bar u$ is approximated by a piecewise constant 
%function, with jumps at the points $x_k$.
\item  At each time $t_j$,  $j\geq 0$, the approximate solution $u(t_j,\cdot)$ is piecewise constant with jumps at the points $x_k$:
$$u(t_j,x)~=~u_{j,k}\qquad\qquad \hbox{for}\quad x_k\leq x< x_{j,k+1}\,.$$
\item For $t_j \leq t < t_{j+1}$ the solution is computed by solving the corresponding
Riemann problems at each point of jump $(t_j,x_k)$, for every integer $k$.
\item  At time $t_{j+1}$ the solution is again approximated by a piecewise constant
function, and the procedure can repeat.
\endi

\begin{figure}[htbp]
\centering
  \includegraphics[scale=0.35]{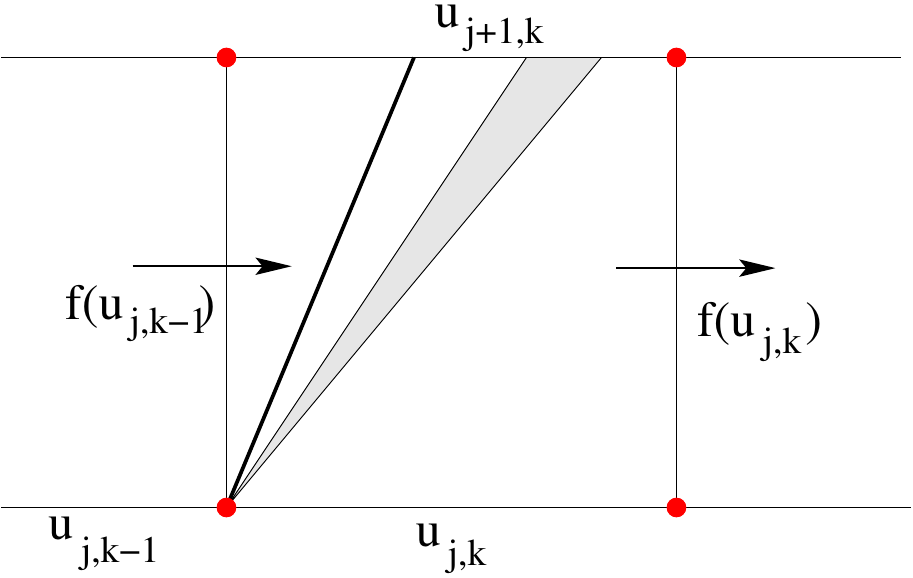}
    \caption{\small By the conservation law (\ref{hcl}), the average value 
    of the solution $u$ 
  on the top side of the square is equal to the value $u_{j+1,k}$
    computed by the formula (\ref{gods}).
    }
\label{f:claw52}
\end{figure}

In the Godunov scheme, at the time $t_{j+1}$ the  function 
$$u(t_{j+1}-, \cdot)~\doteq~\lim_{t\to t_{j+1}-} u(t,\cdot)$$
is replaced by a piecewise constant function,  equal to 
its average on each interval $[x_k, x_{k+1}]$.  Namely
\bel{av}u_{j+1,k}~\doteq~{1\over x_{k+1}- x_k} \int_{x_k}^{x_{k+1}} u(t_{j+1}-, \,
x)\, dx\,.\eeq
A remarkable property of this scheme is that, in order to compute $u_{{j+1},k}$
there is no need to actually construct the solution to a Riemann problem.
Indeed, applying  the divergence theorem on each square of the grid (see Fig.~\ref{f:claw52}),
by the conservation law (\ref{hcl}) one immediately obtains
\bel{gods}u_{j+1,k}~=~u_{j,k} + f(u_{j,k-1}) - f(u_{j, k}).\eeq
The finite difference scheme (\ref{gods}) is called  the (upwind) Godunov scheme.
While this is easy to implement numerically, a rigorous convergence analysis of this scheme
is still lacking, for general solutions containing shocks.

\v
{\bf 2. The Glimm scheme.} 
This scheme is similar to the Godunov scheme, with one major difference.
At each time $t_j$, we have to replace the function $u(t_j-,\cdot)$
with a piecewise constant function having jumps at the 
points $x_k$.    In the Godunov scheme this is achieved by computing the average value
 on each interval $[x_k, x_{k+1}]$.   
In the Glimm scheme~\cite{Glimm}, the restarting is achieved by random sampling. 
Namely, inside each interval $[x_{k}, x_{k+1}]$ we choose a random point $x_{j,k}^*$.
The value $u(t_j-, x_{j,k}^*)$ of the solution to the Riemann problem at this particular point 
is taken to be the value  of the function  $u(t_j, \cdot)$ 
on the entire interval.

More precisely: 
\begi\item 
We choose a random sequence of numbers $\theta_1,\theta_2,\theta_3,\ldots$
uniformly distributed on $[0,1]$.

\item At each time $t_j$, for every $k$ we consider the random point 
$x_{j,k}^* = (k + \theta_j ) \cdot \Delta x$ and define
\bel{restart}u(t_j, x)~=~u(t_j-,  x_{j,k}^*)  \qquad\forall x\in [x_k, \, x_{k+1}[\,.\eeq
\endi

\begin{figure}[htbp]
\centering
  \includegraphics[scale=0.4]{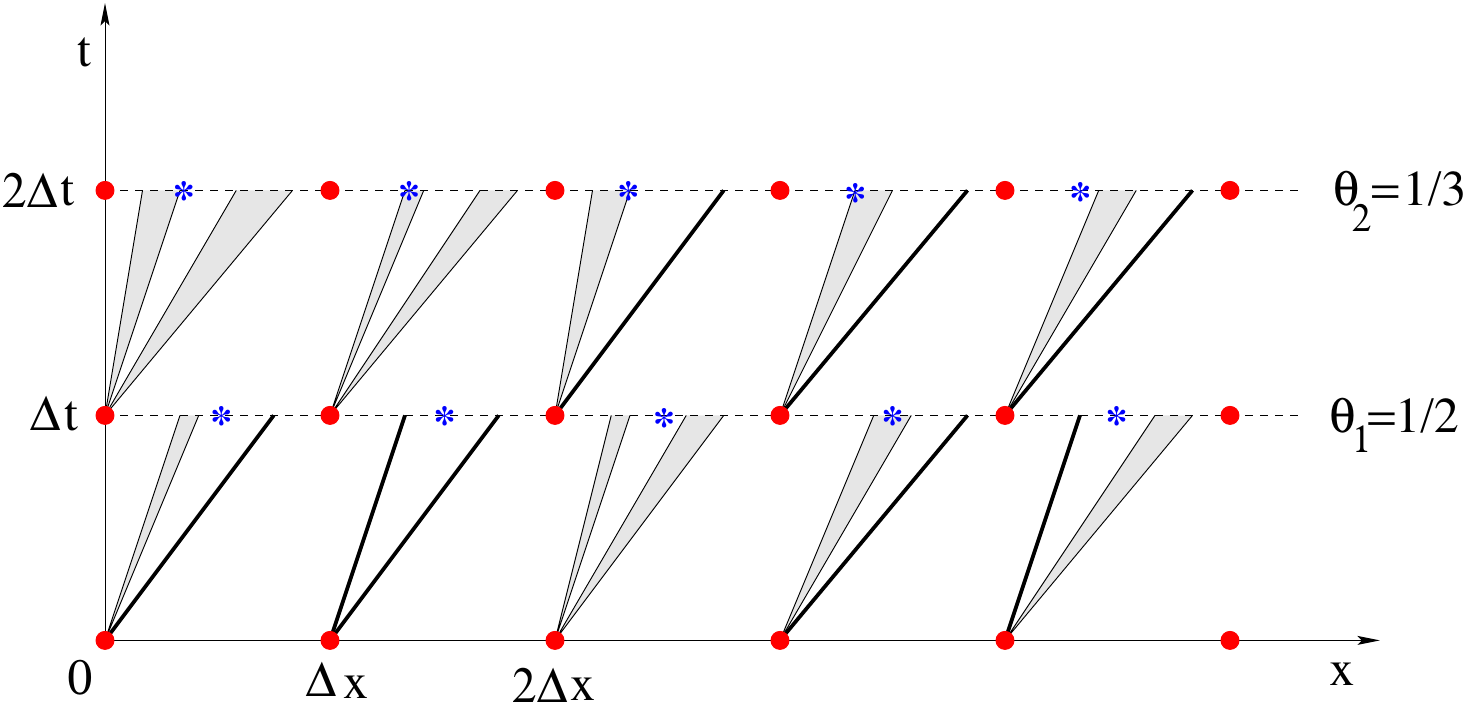}
    \caption{\small An approximate solution obtained by solving a  Riemann problem at each node of the grid. In the Glimm scheme, at each time $t_j = j\, \Delta t$ the solution is sampled at the points marked by an asterisk, depending om the random sequence $\theta_1, \theta_2,\ldots$  
    }
\label{f:claw47}
\end{figure}

As later proved by T.P.Liu \cite{Liu75}, instead of a random sequence one can use 
 a  deterministic sequence of numbers
$\theta_1,\theta_2,\theta_3,\ldots \in [0,1]$  which is    {\bf uniformly distributed},
so that
\bel{unif}\lim_{N\to\infty}
{ \#\big\{ j~;~~1\leq j\leq N,~~\theta_j\in [0,\lambda]~\big\}
\over N}~=~\lambda\qquad\quad \hbox{for each}~\lambda\in [0,1].
\eeq
\begin{figure}[htbp]
\centering
  \includegraphics[scale=0.4]{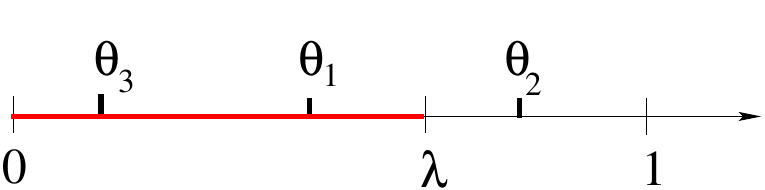}
    \caption{\small A sequence of numbers in the interval $[0,1]$. For every $\lambda$, the percentage of points $\theta_i$, 
    $1\leq i\leq N$ 
    which fall inside the subinterval $[0,\lambda]$ should approach 
    $\lambda$ as $N\to \infty$.}
\label{f:f137}
\end{figure}

A simple way of generating such a sequence is to write decimal digits in inverse order:
\bel{rap}
\theta_1=0.1\,,~~~\ldots~~~,~~~ \theta_{759}=0.957\,,~~~\ldots~~~,~~~
\theta_{39022}=
0.22093\,,~~~\ldots\eeq
\v

\begin{figure}[htbp]
\centering
  \includegraphics[scale=0.4]{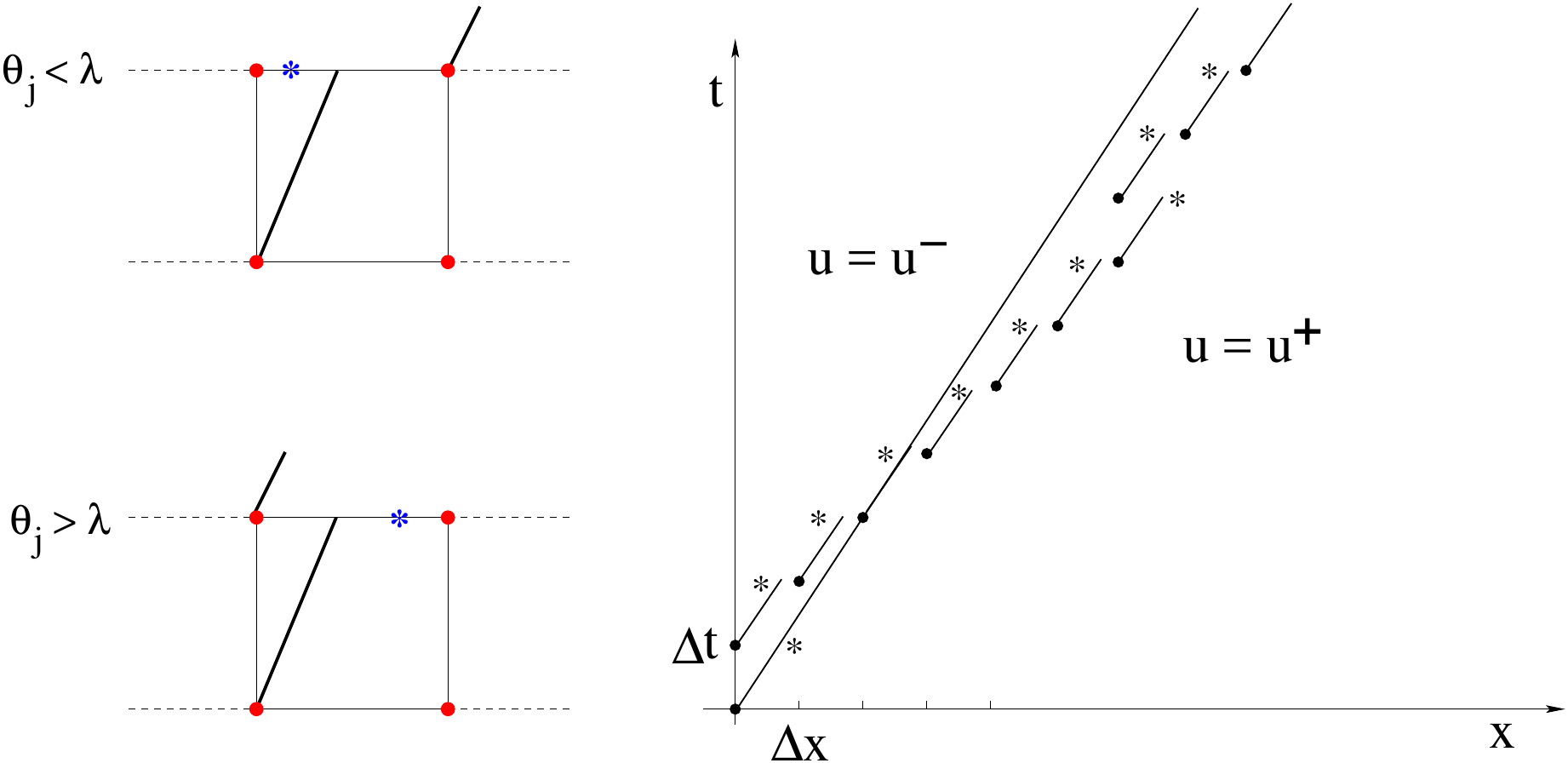}
    \caption{\small Left: at each time step, the position of the shock remains unchanged if $\theta_j>\lambda$, while it jumps forward if $\theta_j<\lambda$.  
    Right: a single shock solution with speed $\lambda\in [0,1]$ and an approximate solution constructed by the Glimm scheme.  As the grid size approaches zero, convergence  to the exact solution is achieved if and only if the limit in (\ref{unif}) holds.}
\label{f:f115}
\end{figure}

The relevance of the assumption (\ref{unif}) is illustrated in Fig.~\ref{f:f115}. 
As in (\ref{shock}), consider a solution containing a single shock, traveling with speed 
$\lambda\in\,]0,1[\,$.

Fix a time interval $[0,T]$ and take $\Delta x = \Delta t = T/N$.   Call $x(t_j)$ the location of the shock at time $t_j$, in the approximate solution.  By construction (see Fig.~\ref{f:f115}, left), we have 
\begi\item If $\theta_j\leq \lambda$, then $x(t_j) = x(t_{j-1}) + \Delta x$.
\item If $\theta_j> \lambda$, then $x(t_j) = x(t_{j-1}) $.
\endi
Hence, at time $T$, the position of the shock in the approximate solution is
$$\bega{rl}
x(T)&=~\#\big\{ j~;~~1\leq j\leq N,~~\theta_j\in
[0,\lambda]~\big\}\cdot
\Delta t \cr&\cr
&=~\displaystyle
{\#\big\{ j~;~~1\leq j\leq N,~~\theta_j\in [0,\lambda]~\big\}
\over N}\cdot T~\to~ \lambda T\qquad\qquad \hbox{as} ~~N\to\infty,
\enda$$
provided that the uniform distribution assumption (\ref{unif}) holds.

\v
{\bf 3. Front tracking approximations.}
In the Glimm scheme, the Riemann problems 
are solved on a fixed grid
in the $t$-$x$ plane.  In a front tracking algorithm, 
the points where new Riemann problems are solved depend on the solution itself. 

\begin{figure}[htbp]
\centering
  \includegraphics[scale=0.4]{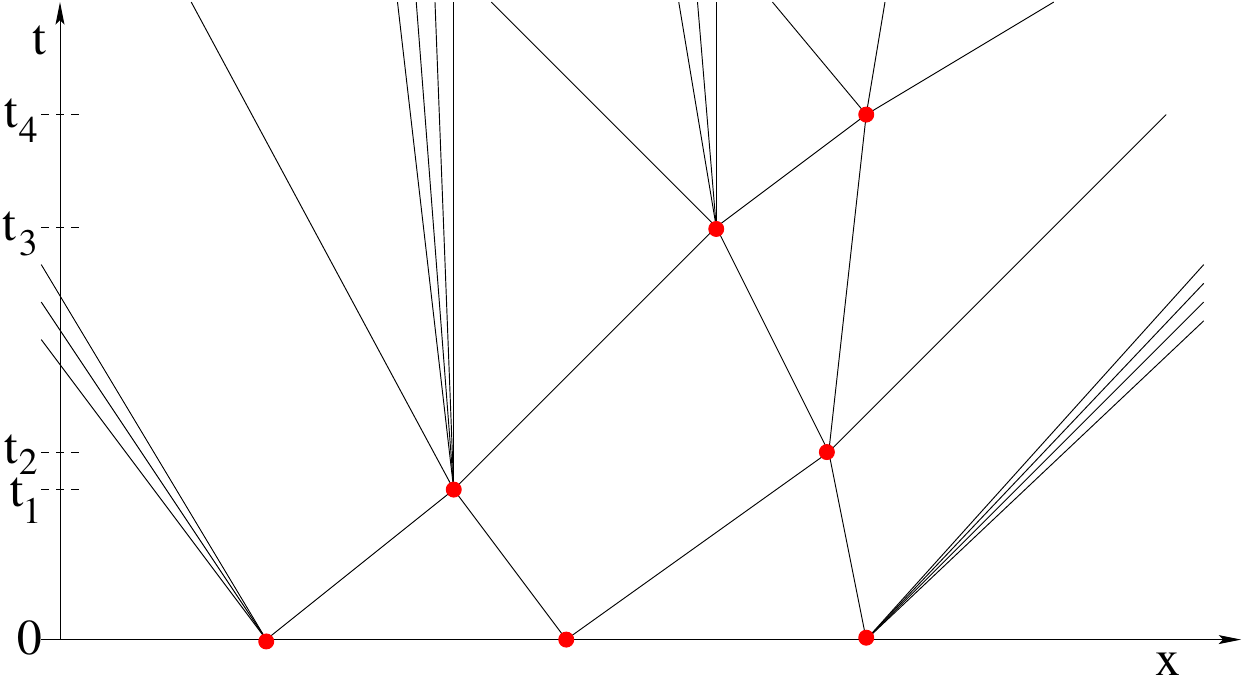}
    \caption{\small An approximate solution constructed by front tracking.
%    Centered rarefaction waves are approximated by piecewise constant functions.
%    New Riemann problems generated by the interaction of two fronts are 
%   approximately solved within a class of piecewise constant functions.
}
\label{f:claw42}
\end{figure}

\begi
\item The construction starts by approximating the initial data $\bar u\in \L^1$ with a piecewise constant function $\bar v$.
\item At each point where $\bar v$ has a jump, an approximate solution to the Riemann problem is constructed, within the class of piecewise constant functions.
As shown in Fig.~\ref{f:claw42}, this  solution can be prolonged 
up to the first time  $t_1$ where two fronts interact.
\item At time $t_1$, we
construct a piecewise constant approximate solution to the new Riemann problem
 generated by the interaction, and prolong the solution until a further interaction occurs.
\item By inductively solving the new Riemann problems at the times $t_2, t_3,\ldots$ where two fronts interact, a piecewise constant approximate solution is constructed for all times $t\geq 0$.
\endi

For a single conservation law, the front tracking method was first  introduced by 
Dafermos~\cite{D72}. 
To apply this method to $n\times n$ systems, one needs to make sure that 
the number of wave fronts does not become infinite in finite time. 
This requires some technical provision, such as the introduction
of ``non-physical fronts" \cite{BJ, Bft}.  
For a comprehensive presentation we refer to \cite{Bbook, Dbook, HR}.

%
%replace centered rarefaction waves
%with piecewise constant rarefaction fans
%
%\begin{figure}[htbp]
%\centering
%  \includegraphics[scale=0.4]{k13.pdf}
%    \caption{\small }
%\label{f:k13}
%\end{figure}
%
%\begi 
%\item The total variation of the approximate solution remains small.
%\item The total number of wave fronts does not become infinite .
%\endi

\v
{\bf 4. Vanishing viscosity approximations.}
Starting with the hyperbolic system (\ref{hcl}) and adding a 
small diffusion term, one obtains the quasilinear parabolic system:
\bel{diffa} 
u^\ve_t + A(u^\ve) u^\ve_x~=~\ve u^\ve_{xx},\qquad\qquad u(0,x) = \bar u(x).\eeq
Here $A(u) = Df(u)$ is the $n\times n$ Jacobian matrix of the flux function.
Letting $\ve\to 0+$, it is expected that the solutions to (\ref{diffa}) will converge to the
unique solution to the hyperbolic Cauchy preblem (\ref{CP3}).   For initial data with small total variation, a rigorous proof 
of this convergence was given in \cite{BiB}.
\v

{\bf 5. Jin-Xin relaxation approximations.} 
These are obtained by solving the second order wave equation
$$u_t + f(u)_x~=~\ve (u_{xx}- u_{tt}) \,,\qquad\qquad u(0,x)= \bar u(x)$$

As $\ve\to 0$, uniform BV bounds and convergence to a unique limit
have been proved by S.~Bianchini in \cite{Bi06}.
The convergence rate has not been studied in detail.

\v
{\bf 6. The method of lines.}
In this case, approximate solutions are obtained by 
discretizing space while keeping time continuous.

\begin{figure}[htbp]
\centering
  \includegraphics[scale=0.4]{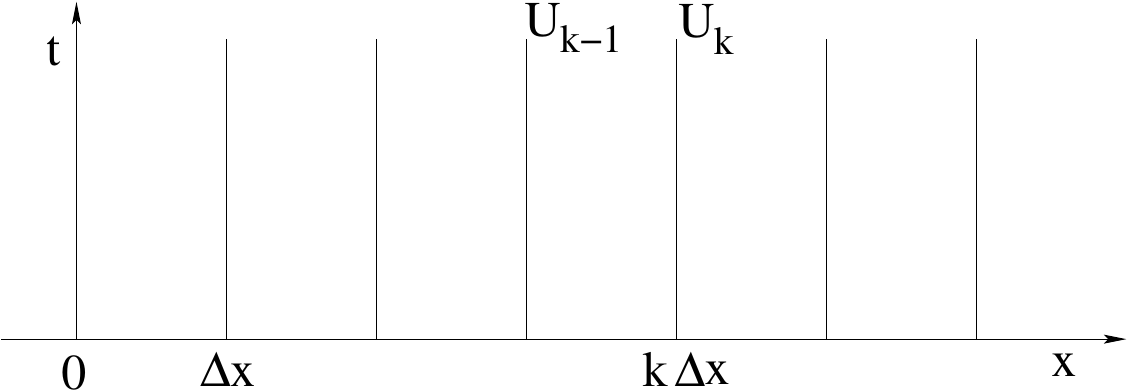}
    \caption{\small The method of lines.}
\label{f:b555}
\end{figure}

Fix a mesh size $\Delta x=\ve>0$.   Approximating the partial derivative $f(u)_x$ 
by a finite difference,  the system of conservation laws (\ref{hcl}) 
is replaced by a countable family of ODEs
$$
{d\over dt} U_k(t)~=~{f(U_{k-1}(t))-f(U_k(t))\over\ve}\,,\qquad\qquad k\in {\mathbb Z},$$
 for the variables $U_k(t)\approx u(t,k\ve)$, see Fig.~\ref{f:b555}.

As $\ve\to 0$, uniform BV bounds and convergence to a unique limit
have been proved by S.~Bianchini in \cite{Bi03}.

%%%%%
{\bf 7. Backward Euler approximations.}
In this case we
discretize time  while keeping space continuous.  Choosing 
$ \Delta t=\ve$ as time step, the Backward Euler approximation takes the form
\bel{BackE}u(t+\ve) ~=~u(t) - \ve f(u(t+\ve))_x\,.\eeq
Setting
$$ v(x)\,=\,u(t,x),\qquad w(x)\,=\,u(t+\ve,x),\qquad A(w)=Df(w),$$
at every time step one needs to solve $$w(x)~=~v(x) -\ve Df(w(x))\, w_x(x),$$
which leads to the ODE
\bel{EODE}w'(x)~=~A^{-1}(w(x)) \left({v(x)-w(x)\over \ve}\right).\eeq
By performing a change of coordinates similar to (\ref{cchange}), 
one can assume that all characteristic speeds (i.e., all eigenvalues of the Jacobian 
matrix $A(u)= Df(u)$) are contained inside 
the interval $[1,2]$. This guarantees that  in (\ref{EODE}) all matrices
$A(w) $ have a uniformly  bounded inverse.
For a fixed $\ve>0$, existence and uniqueness of $\L^1$ solutions to (\ref{EODE}) have been studied in \cite{CZ}, together with traveling wave profiles. 
However, as $\ve\to 0$, uniform BV bounds and convergence 
to a unique limit remain an open question.
\v
{\bf 8. Periodic mollifications.}
These approximations are again constructed by discretizing time.
Fix $\ve>0$ and set $t_k = \ve k$, $k=0,1,2,\ldots$
On each subinterval $[t_k, t_{k+1}[\,,$ the function $u$ is defined to be a classical solution:
$$u_t + Df(u) u_x~=~0,\qquad\qquad t\in [t_k, t_{k+1}[\,.$$
At each time $t_k$, before any shock is formed, the solution is restarted by performing 
a convolution with a mollifying kernel:
$$u(t_k)~=~J_\ve * u(t_k-).$$

Letting $\ve\to 0$, we expect that the approximations should converge to an admissible
weak solution to the original system (\ref{hcl}). This is well known in the scalar case
\cite{HR}. However, 
for general $n\times n$ systems, uniform 
BV bounds and convergence to a unique limit have not been proved.
\v
{\bf 9. Nonlinear  diffusion approximations.}
These take the form
\bel{rvis}u_t + f(u)_x~=~\ve \bigl(B(u)u_x\bigr)_x\,,\eeq
where $B(u)$ is a (possibly degenerate) $n\times n$ diffusion matrix.
Since in many physical systems the viscosity depends on the macroscopic variables,
it would be of great interest to prove rigorous convergence results for solutions of
(\ref{rvis}).   However, apart from the case where $B$ is a constant, invertible matrix
\cite{BiB}, establishing uniform BV bounds and convergence to a unique limit as $\ve\to 0$
remains a challenging open problem.

%%%%%%

\section{Global existence of weak solutions}
\label{sec:4}
\setcounter{equation}{0}
If the initial datum $\bar u:\R\mapsto\R^n$ is smooth, a unique local in time $\C^1$ solution
to the Cauchy problem (\ref{CP3}) can be constructed by the method of characteristic, 
as the fixed point of a contractive transformation \cite{Bbook, RY}.
However, for large times, as shown in Fig.~\ref{f:claw18} the gradient of the solution can 
become unbounded.  Global in time solutions can only be obtained in a space of discontinuous functions.   

For scalar conservation laws, also in several space dimensions, 
a general existence-uniqueness theorem was
proved in the famous paper by Kruzhkov~\cite{Kru}.
Shortly afterwards, an alternative proof based on
the theory of nonlinear contractive semigroups was given by 
Crandall~\cite{Crandall}.

For $n\times n$ hyperbolic systems,
the first global existence theorem for weak solutions to the Cauchy problem (\ref{CP3}) 
was proved in a celebrated paper by Glimm~\cite{Glimm}.

\begin{theorem}\label{t:Glimm} {\bf (Global existence of weak solutions).}
Consider the one-dimensional Cauchy problem (\ref{CP3}) for a  strictly
hyperbolic system of conservation laws, where
each characteristic field is either
linearly degenerate or genuinely nonlinear.

Then there exists a constant
$\delta>0$ such that, if  the initial data 
satisfies 
$$\bar
u\in\L^1(\R;~\R^n),\qquad\qquad \tv \{\bar u\}\leq \delta,$$
then (\ref{CP3}) has a
weak solution $u=u(t,x)$ defined for all $t\geq 0$.

If the system admits a  convex entropy, a global solution exists which is entropy admissible.
\end{theorem}

The proof is achieved by constructing a sequence of approximate solutions $(u_m)_{m\geq 1}$ according to the Glimm scheme (see Fig.~\ref{f:claw47}).   Here we let the grid size
$\Delta t, \Delta x\to 0$  as $m\to\infty$.
\begi
\item
By carefully estimating the strength 
of new waves produced by nonlinear wave interactions, one obtains a uniform bound 
on the total variation of all approximate solutions.  Namely if  $\tv \{\bar u\}$ is sufficiently small, then $\tv \{u_m(t,\cdot)\}$ remains small for all times $t\geq 0$ and every $m\geq 1$.  
\item 
By Helly's compactness theorem, one obtains a convergent subsequence $u_m\to u$
in $\L^1_{loc}(\R_+\times\R\,;~\R^n)$.  
\item
Relying on the assumption that the sequence $(\theta_k)_{k\geq 1}$ used for random 
sampling (\ref{restart})
is uniformly distributed on $[0,1]$, one proves that with probability one the limit function 
$u=u(t,x)$ is a weak solution to the Cauchy problem.
\endi

For many years, the Glimm scheme provided the only tool for 
a rigorous analysis of 
weak solutions to hyperbolic conservation laws.  
Among the first such studies, 
the asymptotic behavior 
of solutions as   $t\to +\infty$ was analyzed by T.P.Liu~\cite{Liu77}.   
The assumption of genuine nonlinearity or linear degeneracy of each characteristic field was removed in \cite{Liu81}.  
Alternative proofs, relying on front tracking approximations, were later 
given in \cite{AM01, BJ, Bft}.

We remark that, in all these results, the smallness of the initial data is a key assumption which has never been removed. This leads to the following question:

{\it If the total variation of the initial data $u(0,\cdot) =\bar u(\cdot)$ is bounded but possibly large,
does the total variation of the solution $u(t,\cdot)$  remain bounded for all times $t>0$, or can it blow up in finite time?}

A counterexample constructed by Jenssen~\cite{J} shows that, in some cases,
the total variation of a weak solution can indeed blow up in finite time.
However, the hyperbolic system considered in this example does 
not admit any strictly convex entropy.  In particular, the construction does not apply
to any of the systems of conservation laws which are 
relevant for continuum physics.

{\bf Open Problem \#1.} {\it
Consider a strictly hyperbolic system
of conservation laws (\ref{hcl}), admitting a strictly convex entropy.
\begi
\item Construct an example  of an entropy admissible weak solution $u=u(t,x)$ 
whose total variation blows up in finite time.

 \item Or else, prove that every solution, whose total variation is initially bounded, remains
 with bounded variation for all times $t>0$.
\endi
}

\begin{figure}[htbp]
\centering
  \includegraphics[scale=0.3]{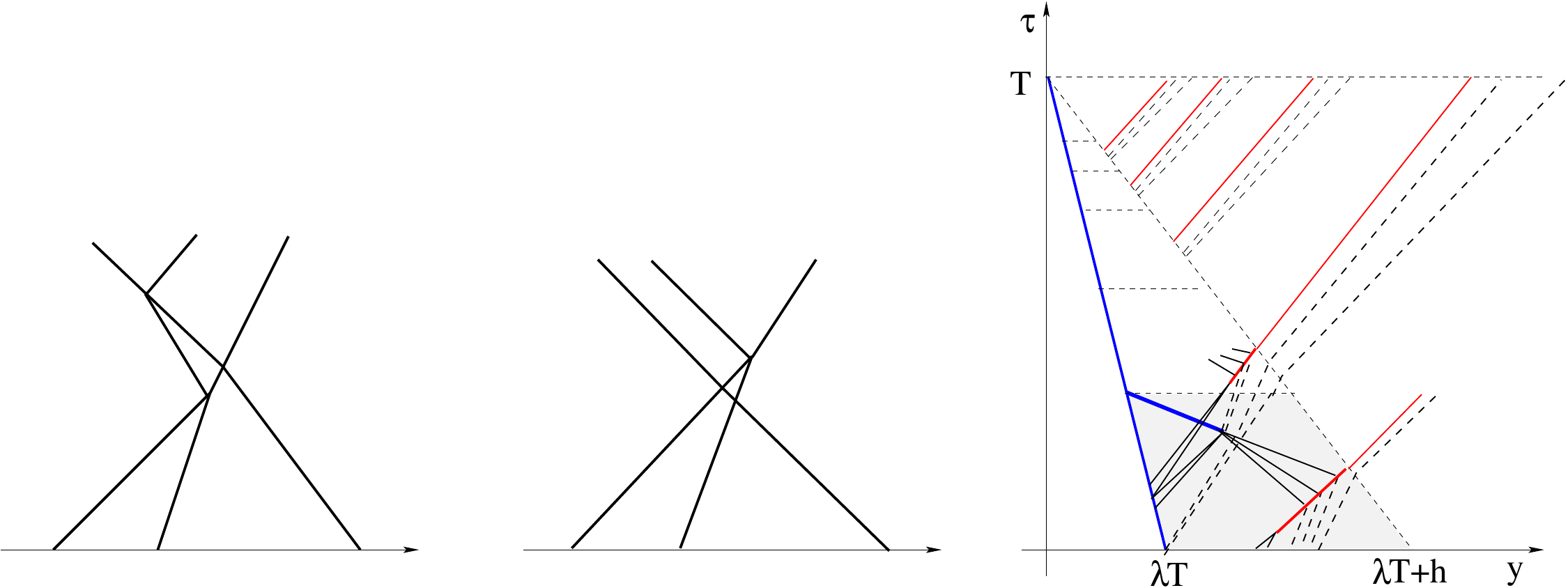}
    \caption{\small Left: if their speeds are slightly changed, the three wave
   fronts  will interact in different order, producing new waves of different strengths.  
    Right: a sketch of the interaction pattern considered in \cite{BCZ}, leading to  blow up of the total variation in finite time.}
\label{f:claw79}
\end{figure}

The question of global BV bounds versus finite time blow up is not resolved even for the
$2\times 2$ system of isentropic gas dynamics (\ref{igd}).  The recent analysis in \cite{BCZ}
only points out how difficult the problem really is.   
By slightly changing the wave speeds,
one can arrange so that wave fronts cross each other in a different order
(see Fig.~\ref{f:claw79}). 
For the same initial data,
one can construct approximate solutions whose total variation remains bounded,
and other approximate solutions whose total variation blows up in finite time, 
depending on the interaction pattern.  
It is hard to say what happens for the exact solution.

\section{Continuous dependence on initial data}
\label{sec:5}
\setcounter{equation}{0}

For a wide class of evolution equations, continuous dependence on initial data
is achieved by showing that the distance between any two solutions 
satisfies the differential inequality 
\bel{gin}
{d\over dt} \bigl\| u(t)-v(t)\bigr\|~\leq~C \bigl\| u(t)-v(t)\bigr\|,\eeq
for some constant $C$.
In turn, by Gronwall's lemma this implies
\bel{Gro2}\bigl\| u(t)-v(t)\bigr\|~\leq~e^{Ct}\, \bigl\| u(0)-v(0)\bigr\|\,.\eeq
%In particular, this is the standard approach in the theory of ODEs
%with Lipschitz continuous right hand side.

This approach can be applied to scalar conservation laws, 
even in a multidimensional space $\R^d$.  As proved in \cite{Crandall, Kru},
a scalar conservation law generates a contractive semigroup on $\L^1(\R^d)$.
Namely,
 for every couple of entropy admissible solutions $u,v$ one has
\bel{contsem} 
\bigl\| u(t)-v(t)\bigr\|_{\L^1}~\leq~ \bigl\| u(s)-v(s)\bigr\|_{\L^1}\qquad\qquad\forall 
0\leq s<t.\eeq

\begin{figure}[htbp]
\centering
  \includegraphics[scale=0.4]{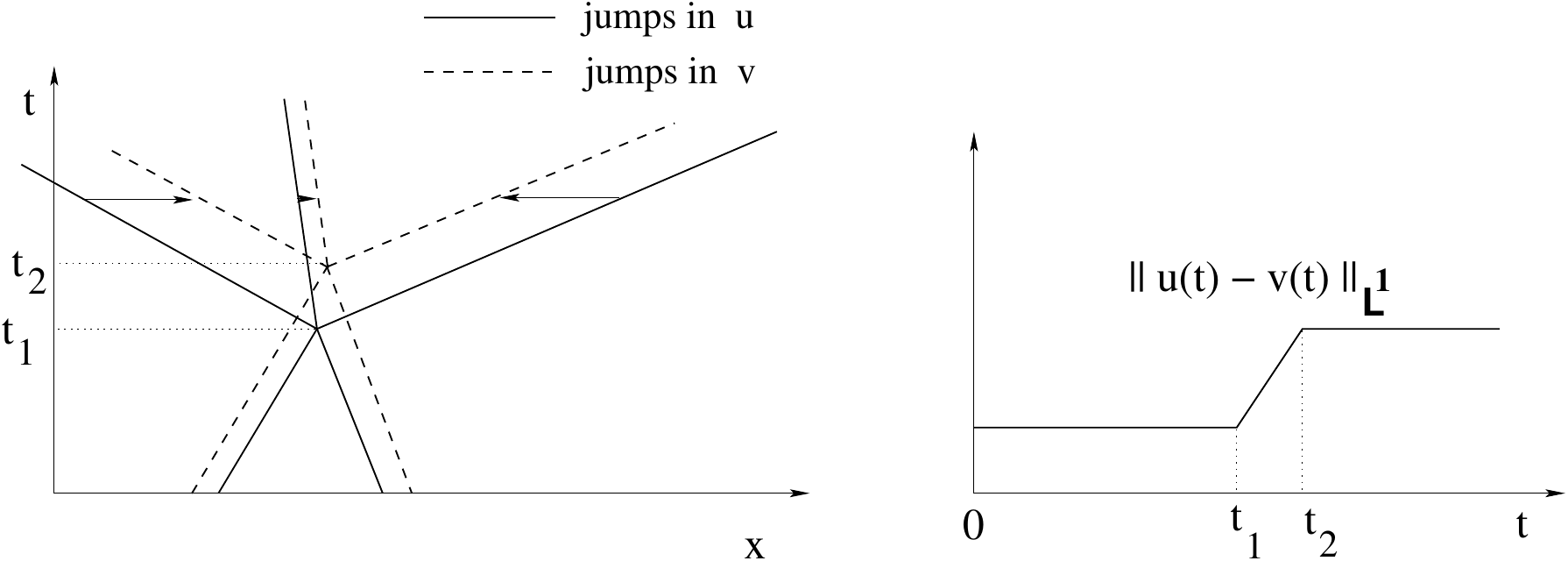}
    \caption{\small Left: a solution $u$ which initially contains two approaching shocks, and a second solution $v$ which differs from $u$ only in the location of the shocks.  
    During the short  time interval $[t_1, t_2]$ the solution $u$ contains
three shocks, while $v$ still has two.   Right: the $\L^1$ distance
between the two solutions remains constant for $t\in [0,t_1]$ and for $t\geq t_2$, but increases rapidly during the interval $[t_1, t_2]$.
    }
\label{f:h67}
\end{figure}

On the other hand, the inequality (\ref{gin}) does not hold for weak
solutions to hyperbolic systems. As shown in Fig.~\ref{f:h67},
the $\L^1$ distance between two nearby solutions can increase
rapidly during short time intervals.
\v

\begin{figure}[htbp]
\centering
  \includegraphics[scale=0.4]{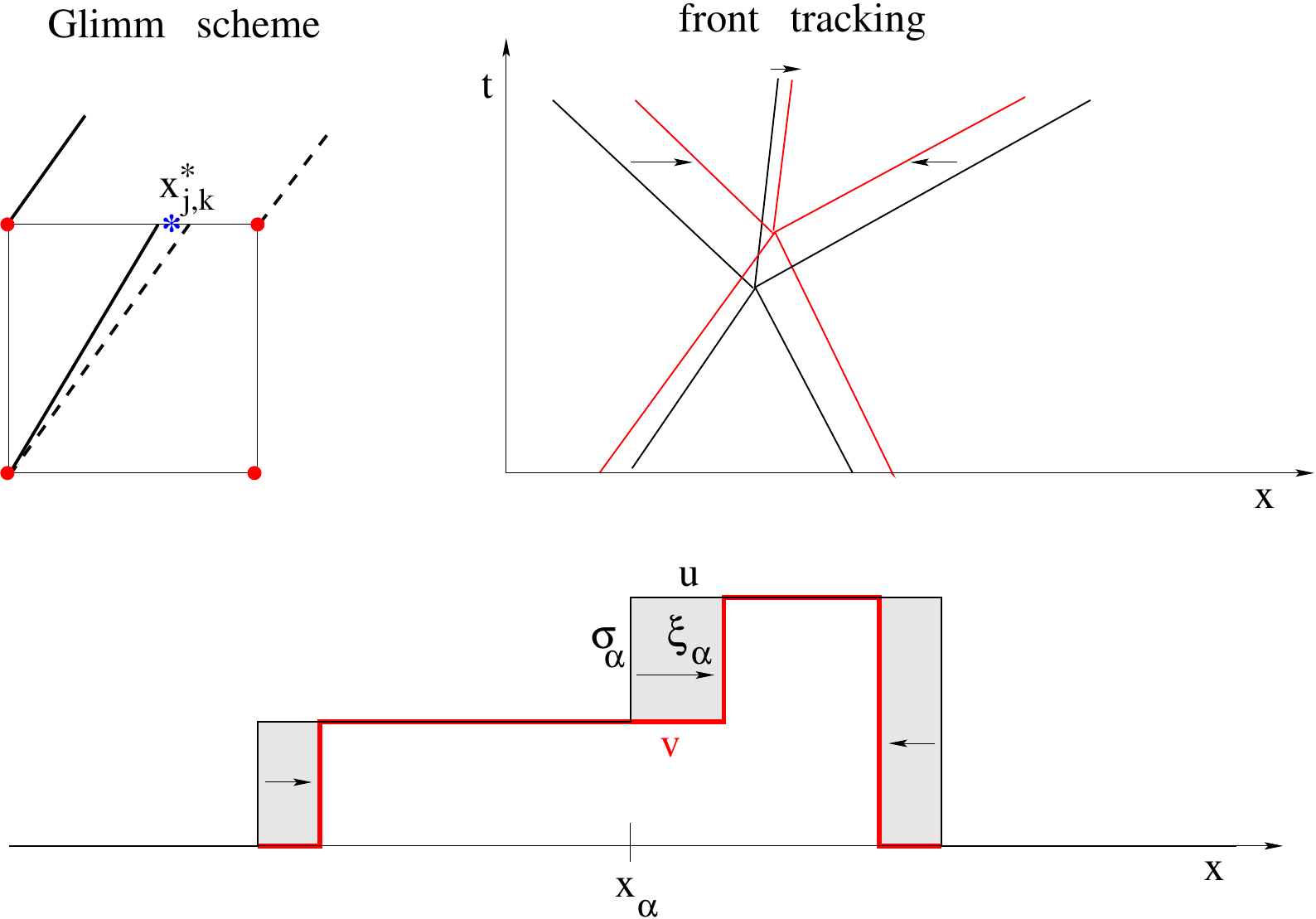}
    \caption{\small For the Glimm scheme, approximate solutions do not depend continuously on the initial data.  On the other hand, by shifting the location of the wave fronts,
    one can understand how front tracking approximations are affected 
    by changes in the initial data.}
\label{f:claw54}
\end{figure}

We observe that the Glimm scheme does not provide insight on the continuous
dependence of solutions.
Indeed, the approximate solutions constructed by the Glimm scheme do not
depend continuously on initial data. With reference to Fig.~\ref{f:claw54}, left,
an arbitrarily small change in the speed of a shock may place it to the right or to the left of the sampling point $x^*_{j,k}$.   
In this case, the piecewise constant approximation changes value on the entire 
interval $[x_k, x_{k+1}]$.

Understanding how weak solutions depend on initial data was the primary motivation for developing an alternative approximation scheme based on wave front tracking
 \cite{Bft}.  As shown in Fig.~\ref{f:claw54}, let $u$ be a piecewise constant 
 approximate solution, and let $v$ be a perturbed solution obtained by slightly shifting 
 the location of the jumps in $u$.  At any time $t\geq 0$, the $\L^1$ distance between 
 the two solutions is measured by
 \bel{jss}\bigl\|u(t)-v(t)\bigr\|_{\L^1}~\approx~\sum_\alpha |\sigma_\alpha|\cdot |\xi_\alpha|~
=~\sum_\alpha \hbox{[jump strength]}\times\hbox{[shift]}.\eeq
By carefully estimating how the right hand side of (\ref{jss}) changes at 
interaction times, one obtains a bound on the distance between the two solutions,
at every time $t\geq 0$.

\begin{theorem} Let  (\ref{hcl}) be a strictly hyperbolic $n\times n$ system of conservation laws, and assume that each 
characteristic field is either genuinely nonlinear or linearly degenerate.  
Then 
there exists a domain $\D\subset \L^1(\R;\,\R^n)$ and a semigroup
$S:\D\times \R_+\mapsto\D$ with the following properties.
\begi
\item[(i)] The domain $\D$ contains all functions $\bar u\in \L^1$ with sufficiently small total variation.
\item[(ii)]
The map
$(\bar u,t)\mapsto u(t,\cdot)\doteq S_t\bar u$ is a uniformly
Lipschitz continuous semigroup: 
$$S_0\bar u=\bar u,\qquad S_s (S_t\bar u)=S_{s+t}\bar
u,$$
$$\big\|S_t\bar u-S_s\bar v\big\|_{\L^1}~\leq~ L\cdot \big(\|\bar u-
\bar v\|_{\L^1}+|t-s|\big)\qquad\quad\hbox{for all}~~\bar u,\bar
v\in\D,~~s,t\geq 0.$$
\item[(iii)] For every initial data 
$\bar u\in \D$, the  trajectory $t\mapsto u(t,\cdot) = S_t\bar u$
provides an admissible solution to the Cauchy problem (\ref{CP3}).
\endi
\end{theorem}

Trajectories of the semigroup can be obtained as limits of front tracking approximations.
This theorem was first proved in \cite{BC1} in the case of $2\times 2$ systems.
A more elaborate proof, valid for $n\times n$ systems, was later worked out in \cite{BCP}.
These earlier proofs relied on a homotopy method: the distance between two solutions
$u,v$ was estimated by constructing a 1-parameter family of solutions $u^\theta$, $\theta\in [0,1]$, connecting $u$ with $v$. At every time $t\geq 0$, the  
length of the path $\theta\mapsto u^\theta(t,\cdot)$ provides
a bound on the distance $\bigl\|u(t)-v(t)\bigr\|_{\L^1}$.

\begin{figure}[htbp]
\centering
  \includegraphics[scale=0.4]{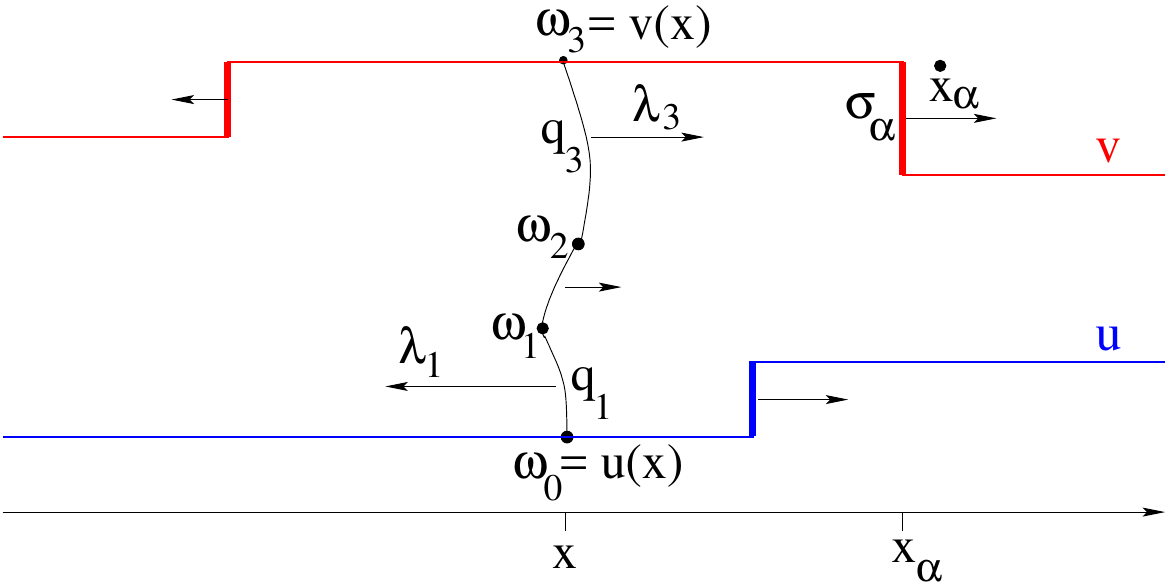}
    \caption{\small  Constructing the functional $\Phi(u,v)$.
    At each point $x$, the strengths $q_1(x), \ldots, q_n(x)$ of the shocks connecting $u(x)$ with $v(x)$ can be regarded as the  scalar components of the jump $\bigl(u(x), v(x)\bigr)$.}
\label{f:p9}
\end{figure}

Using ideas introduced by T.P.Liu and T.Yang \cite{LY},
an alternative proof was given in \cite{BLY}.  
This approach relies 
on the construction of a Lyapunov functional $\Phi: \D\times\D\mapsto \R_+$ with the following properties.
\begi
\item~~ $\Phi$ is equivalent to the $\L^1$ distance:
$${1\over C}\cdot
\big\|v-u\big\|_{\L^1~}\leq~ \Phi(u,v)~\leq~ C\cdot
\big\|v-u\big\|_{\L^1}\,,$$
for every couple of piecewise constant functions $u,v\in \L^1(\R;\,\R^n)$ with small total variation.
\item ~~$\Phi$ is non-increasing in time, along couples of (front tracking approximate) solutions:
\bel{Pdec} \Phi(u(t), v(t))~\leq~ \Phi(u(s), v(s))\qquad\qquad \forall ~t\geq s\geq 0.\eeq
\endi
Given two piecewise constant functions $u,v$, the functional $\Phi(u,v)$ is defined
as follows.
For each $x\in\R$, we uniquely determine intermediate states
$$u(x)=\omega_0(x),\quad \omega_1(x),\quad\ldots,\quad\omega_n(x)=v(x),$$
such that every pair $\bigl(\omega_{i-1}(x),\, \omega_i(x)\bigr)$ is joined by a 
(possibly non-admissible)
shock of the $i$-th family (see Fig.~\ref{f:p9}).   Calling $q_i(x)$ the strength of this shock,
the $\L^1$ distance between $u$ and $v$ can now be estimated  as 
$$\|u-v\|_{\L^1}~\approx~\int_{-\infty}^{+\infty} \sum_{i=1}^n |q_i(x)|\, dx\,.$$
To achieve the decreasing property (\ref{Pdec}), suitable weights $W_i$ must be inserted.  
Roughly speaking, $W_i(x)$ measures the total strength of waves in $u$ and in $v$ that approach an $i$-shock located at $x$.  The functional $\Phi$ thus takes the form
$$\Phi(u,v)~=~\sum_{i=1}^n \int_{-\infty}^\infty
|q_i(x)|\, W_i(x)\, dx\,.$$
For all details we refer to \cite{BLY}.   This approach greatly simplified the earlier proofs, 
and is now adopted in most textbooks on the subject \cite{Bbook, Dbook, HR}.

\v

All previous results were proved in the same setting as in Glimm's theorem, where
each characteristic field is either linearly degenerate or genuinely nonlinear.
Eventually, this assumption was entirely removed by the approach based on 
vanishing viscosity
approximations \cite{BiB}. In this case, one does not even need that 
the hyperbolic system be in conservation form.

\begin{theorem} \label{t:52}
Consider the Cauchy problem for a strictly hyperbolic
system with small viscosity
\bel{evisc}u_t+A(u)u_x=\ve\,u_{xx}\,,\qquad\qquad u(0,x)=\bar u(x)\,.\eeq
If
$\tv\{\bar u\}$ is sufficiently small,
then (\ref{evisc}) admits a unique solution~~
$u^\ve(t,\cdot)=S^\ve_t\bar u$, defined for all $t\geq 0$.
Moreover, for some constants $C,L$ independent of $\ve$, one has 
$$\tv\big\{S_t^\ve \bar u\big\}~\leq~ C\,
\tv\{\bar u\}\,,\eqno{\bf (BV~bounds)}$$
$$\big\|S^\ve_t\bar u-S^\ve_t\bar v\big\|_{\L^1}~\leq~ L\,
\|\bar u-\bar v\|_{\L^1}\,.\quad ~\eqno{\bf  (\L^1 ~stability)}$$
As $\ve\to 0$, the solutions $u^\ve$ converge to the trajectories of
a semigroup $S$ such that
$$\big\|S_t\bar u-S_t\bar v\big\|_{\L^1}~\leq ~L\, \|\bar u-\bar v\|_{\L^1}
\qquad\qquad \forall t\geq 0\,.$$

If the system is in conservation form: $A(u)=Df(u)$ for some flux function $f$, 
then every trajectory $t\mapsto u(t,\cdot) = S_t\bar u$ provides a weak solution to the 
Cauchy problem
$$u_t+f(u)_x~=~0,\qquad\qquad u(0,x)=\bar u(x).$$
Moreover, the Liu  admissibility conditions (\ref{liuadm}) are satisfied at every point 
of approximate jump.
\end{theorem}

For a proof, see \cite{BiB} or the lecture notes \cite{Btut}.
These vanishing viscosity limits
can be regarded as the unique {\bf viscosity solutions}
of the hyperbolic Cauchy problem
$$u_t+A(u)u_x=0\qquad\qquad u(0,x)=\bar u(x)\,.$$
We remark that all  results stated in Theorem~\ref{t:52} hold for the system (\ref{evisc})
with ``artificial viscosity", where the diffusion acts uniformly on all components
$(u_1,u_2,\ldots, u_n)$ of the solution.   In several  physical models,
the diffusion acts differently on different components, and depends on the state $u$
as well.   This leads to

{\bf Open Problem \#2.} {\it  Extend the results of Theorem~\ref{t:52} to 
hyperbolic systems with nonlinear diffusion:
\bel{hvis}u_t + A(u)u_x~=~\ve \bigl(B(u)u_x\bigr)_x\,,\qquad\qquad u(0,\cdot)=\bar u\,,\eeq
where $B(u)$ is a positive semidefinite viscosity matrix.
Assuming that the initial data $\bar u$ has small total variation, prove uniform bounds
on $\tv\{ u(t,\cdot)\bigr\}$ and study the limit of these solutions as $\ve\to 0+$.
}

In the conservative case where $A(u) = Df(u)$, we expect that the vanishing viscosity
limit should be unique, and provide a weak solution to the hyperbolic system (\ref{hcl}).
On the other hand, in the non-conservative case, as $\ve\to 0+$ different limits may well be obtained, depending on the choice of the viscosity matrices $B(u)$ in (\ref{hvis}).

\section{Uniqueness of weak solutions}
\label{sec:6}
\setcounter{equation}{0}

Having constructed a Lipschitz semigroup of admissible weak solutions, which are limits of vanishing viscosity approximations (and of front tracking approximations as well),
it becomes entirely clear which is the unique ``good" solution to the Cauchy problem
(\ref{CP3}).   Namely, the semigroup trajectory $t\mapsto u(t,\cdot)=S_t\bar u$.
From this point of view, uniqueness becomes a marginal issue in the overall theory.
Some authors barely mention the problem \cite{Liubook}, focusing instead all the attention on the continuous dependence on initial data.

On the other hand, a general uniqueness theorem can be quite useful
if we want to study the convergence of different approximation methods. 
Without a uniqueness result,
one may even suspect that these algorithms converge to different limit solutions.

\begin{figure}[htbp]
\centering
  \includegraphics[scale=0.4]{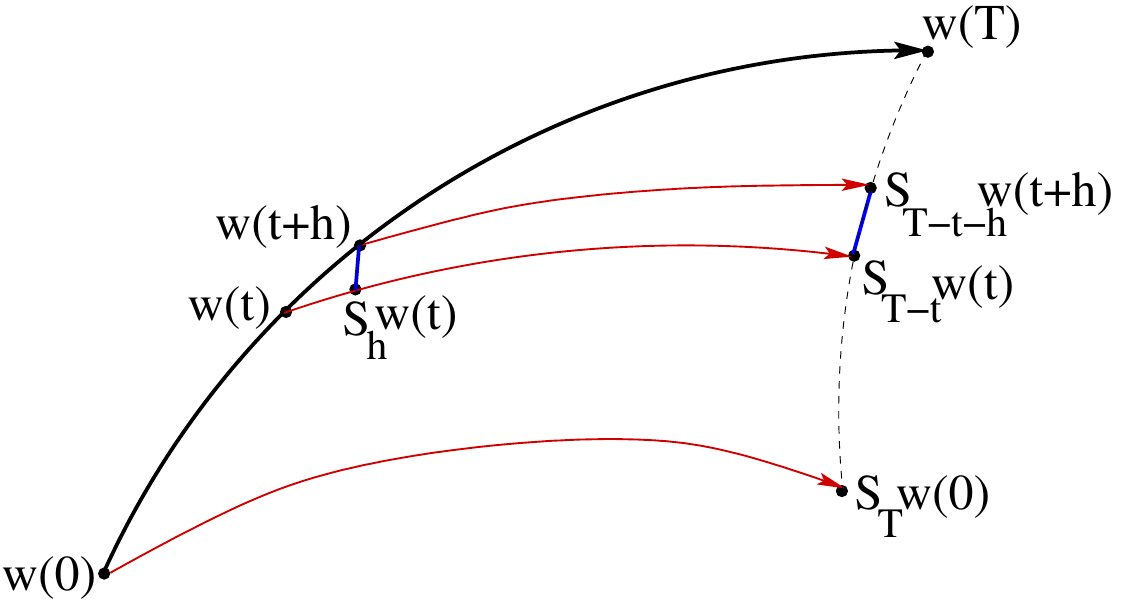}
    \caption{\small The distance between $w(T)$ and $S_T w(0)$ is 
    bounded by the length of the path $t\mapsto S_{T-t}w(t)$, $t\in [0,T]$. Since the semigroup 
    amplifies distances at most by a factor $L$, this leads to the formula (\ref{eest}).}
\label{f:f1}
\end{figure}

As soon as a semigroup of solutions has been constructed,
to establish a uniqueness theorem it suffices to come up with a set of conditions that
uniquely characterize the semigroup trajectories.  Various ways to do this 
have been worked out in \cite{BG1, BLF,  BL}.
The proofs rely on the following elementary error estimate (see Fig.~\ref{f:f1}).
\begin{lemma}
 Let
$S:\D\times [0,\infty[\,\mapsto\D$ be a Lipschitz semigroup
satisfying
 $$\|S_t u-S_s v\|~\leq~ L\cdot\|u-v\|+L'\cdot|t-s|\,.$$

Then, for every Lipschitz continuous map $w:[0,T]\mapsto \D$ one has
\bel{eest}
\big\| w(T)-S_{\strut T} w(0)\big\|~\leq~\displaystyle L\cdot\int_0^T
\left\{\liminf_{h\to 0+}{\big\| w(t+h)-S_hw(t)\big\|
\over h}\right\}\,dt\,.\eeq
\end{lemma}
For a proof, see \cite{BGlimm, Bbook}.
The integrand on the right hand side of (\ref{eest}) can be interpreted as the
{\bf instantaneous error rate} of the approximate solution $w(\cdot)$  at time $t$.

To prove that a solution $u=u(t,x)$ of the Cauchy problem (\ref{CP3}) coincides with the
semigroup trajectory, it now suffices to show that
\bel{linf}
\liminf_{h\to 0+}~{\big\| u(\tau+h)-S_hu(\tau)\big\|_{\L^1}
\over h}~=~0\eeq
for a.e.~time $\tau\in [0,T]$.

\begin{figure}[htbp]
\centering
  \includegraphics[scale=0.4]{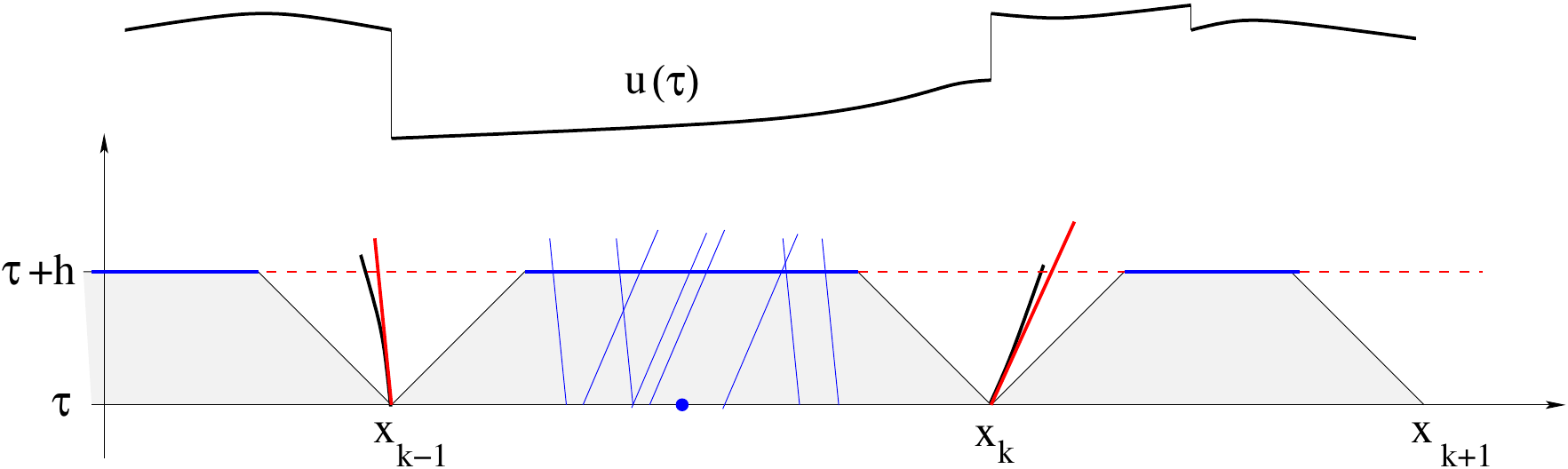}
    \caption{\small Given the function $u(\tau,\cdot)$, we insert points $x_k$ so that
    the total variation on each open interval ~$]x_{k-1}, x_k[\,$ is $<\ve$.}
\label{f:claw81}
\end{figure}

To fix ideas,  w.l.o.g.~we assume that all wave speeds (i.e., all eigenvalues of the matrices $A(u)=Df(u)$) are contained in the interval
$[-1,1]$.    Given the function $u(\tau,\cdot)$, following the approach introduced in \cite{BGlimm} we split the real line inserting 
points 
$$-\infty~\doteq~x_0<x_1<x_2<\cdots<x_N\doteq +\infty,$$ 
so that (see Fig.~\ref{f:claw81})
\bel{tve}\tv\bigl\{u(\tau, \cdot)\,;~]x_{i-1},x_i[\,\bigr\}~<~\ve,\qquad\qquad i=1,2,\ldots,N.\eeq

We now estimate
\bel{er1}\bega{l}\ds
{1\over h}\int_{-\infty}^\infty \bigg|u(\tau+h,x) - S_h u(\tau)(x)\bigg|\, dx\\[4mm]
 \ds =~ \sum_k  {1\over h} \int_{x_k-h}^{x_k+ h}
 \bigg|u(\tau+h,x) - S_h u(\tau)(x)\bigg| \, dx + \sum_k  {1\over h} \int_{x_{k-1} + h}^{x_k- h}
 \bigg|u(\tau+h,x) - S_h u(\tau)(x)\bigg| \, dx\\[4mm]
=~\ds  \sum_k  A_k(h)+  \sum_k B_k(h)\,.\enda\eeq
We claim that, by choosing $\ve>0$ small, the limit as $h\to 0$ of the 
right hand side of (\ref{er1})
can be made arbitrarily small.  Indeed, assume that at $(\tau, x_k)$ the solution $u$ 
is either approximately continuous, or has an approximate jump, as in (\ref{ajump}).
Then
\begi
\item On each interval $[x_k-h, x_k+h]$, the function $u(\tau+h,\cdot)$ and the 
semigroup solution $S_h u(\tau)$ are both compared with the piecewise constant function
$$U_k(t,x)~\doteq~\left\{
\bega{rl}u^+\,\doteq\, u(\tau,x_k+)\qquad &\hbox{if}\qquad x>x_k+\lambda (t-\tau), \cr
u^-\,\doteq\, u(\tau,x_k -)\qquad &\hbox{if}\qquad x<x_k+\lambda (t-\tau). \enda\right.$$
%
%$$\implies\qquad \lim_{h\to 0+}{1\over h}\int_{\xi-h}^{\xi+h}
%\Bigg| u(\tau+h,~x)-U_{(\tau,\xi)}^\sharp(\tau+h,~x)\Bigg|~dx
%~=~0\eqno(E1)$$
\item On each of the remaining intervals $[x_{k-1} + h, x_k - h]$ the function $u(\tau+h,\cdot)$ and the  semigroup solution $S_h u(\tau)$ are both compared with 
the solution $W_k$ of the
linear Cauchy problem with constant coefficients
\bel{linp}
w_t+A_k w_x\,=\,0\qquad\qquad w(\tau,x)=u(\tau,x),\eeq
where $A_k\doteq Df\bigl(u(\tau,y_k)\bigr)$
for some $y_k\in ~]x_{k-1}, x_k[\,$.
\endi

These comparisons are based on two lemmas.  Equivalent results were proved in 
Theorem~2.6 of
 \cite{Bbook} and in \cite{BGu}, respectively.  
 Since they play a crucial role in the uniqueness results,
 we include here the complete proofs.

\begin{lemma} \label{l:62} Assume that the map  $t\mapsto u(t,\cdot)$ is Lipschitz continuous with values in $\L^1(\R)$.   Moreover, let $(\tau,\xi)$ be a point of approximate jump for $u$, so that (\ref{ajump}) holds.
Then
\bel{er2}\lim_{h\to 0+}  {1\over h} \int_{-h}^ h
 \bigg|u(\tau+h,\xi+x) - U(h,x)\bigg| \, dx ~=~0.\eeq
\end{lemma}

{\bf Proof.} Assume that, on the contrary, there exists a decreasing 
sequence $h_m\to 0$ such that
$$ {1\over h_m} \int_{-h_m}^{ h_m}
 \bigg|u(\tau+h_m,\xi+x) - U(h_m,x)\bigg| \, dx ~\geq~\delta~>~0\qquad
 \qquad\forall m\geq 1.$$
By Lipschitz continuity, for some constant $L$ this implies
$$  \int_{-h_m}^{ h_m}
 \bigg|u(\tau+h_m+s,\xi+x) - U(h_m+s,x)\bigg| \, dx ~\geq~\delta h_m-Ls\quad
 \qquad\forall m\geq 1,\, s\geq 0.$$
Setting 
$$\bar s = {\delta\over L}\,, \qquad\qquad r_m\,=\, h_m(1+\bar s),$$
we now obtain
$$\bega{l}\ds {1\over r^2_m} \int_{-r_m}^{r_m} \int_{-r_m}^{r_m}
 \bigg|u(\tau+t,\,\xi+x) - U(t,x)\bigg| \, dx dt\\[4mm]
 \qquad\ds \geq~{1\over (1+\bar s)^2 h^2_m}
 \int_{h_m}^{ (1+\bar s)h_m} \int_{-h_m}^{ h_m}
 \bigg|u(\tau+t,\,\xi+x) - U(t,x)\bigg| \, dx dt\\[4mm]\qquad
 \ds \geq~{1\over (1+\bar s)^2 h^2_m} \int_0^{\bar s h_m}
 (\delta h_m- Ls)\, ds~=~{1\over (1+\bar s)^2} \cdot {\delta^2\over 2 L}~>~0. 
\enda
$$
This contradicts the assumption (\ref{ajump}).
\endproof

%For a proof, see Theorem~2.6 in \cite{Bbook}.

Estimating the difference $u(\tau+h, \cdot) - W_k(\tau+h,\cdot)$ requires more work.
Indeed, here we are approximating $u$ with the solution to a linearized problem.
This is accurate as long as $u$ remains close to a constant, not only at time $\tau$, 
but on the entire trapezoidal domain
\bel{Gade}\Gamma_k~=~\Big\{ (t,x)\,;~~t\in [\tau, \tau+h],~
~x\in  J_k(t)\Big\},\eeq
\bel{Jkdef} J_k(t)\,\doteq~]x_{k-1}+(t-\tau), ~x_k - (t-\tau)\bigr[\,,\eeq
as shown in Fig.~\ref{f:claw81}.
For a solution constructed by front tracking, or by the Glimm scheme,
the assumption that $u(\tau,\cdot)$ has small total variation on $]x_{k-1}, x_k[$
 implies that $u$ has few waves (and hence is nearly constant) also on the 
 domain of dependence $\Gamma_k$ in (\ref{Gade}).  However, in principle 
 this may not be true
 for more general weak solutions.  
 For example, the analysis in \cite{BBJ} shows that the Godunov scheme can ampify
 the total variation by an arbitrarily large factor.
 
  For this reason, additional assumptions
 were required in earlier papers.   Namely ``Tame Variation" in \cite{BLF},
 ``Tame Oscillation" in \cite{BG1}, and ``Bounded variation along space-like curves"
 in \cite{BL}.     
 For a class of $2\times 2$ systems, the recent analysis in \cite{CKV}
 has shown that any entropy admissible weak solution taking values in the domain 
 of the semigroup
 always satisfies  the bounded variation condition in \cite{BL}. 
 Therefore, one does not need this additional assumption to achieve uniqueness.

 Following \cite{BDL, BGu}, we show here that uniqueness holds also  for 
 fully general $n\times n$ systems, without any of the previous 
 regularity assumptions.

To appreciate the underlying idea, consider the  scalar function
 $$V(t)~=~\tv\bigl\{u(t,\cdot)\,;~J_k(t)\bigr\},$$
 where $J_k(t)$ is the interval introduced in (\ref{Jkdef}).
Notice that $V$ is lower semicontinuous, hence measurable. 
Assume that 
$\tau$ is a Lebesgue point for $V$. Since $V(\tau)\leq\ve$, this implies that
the set of times where $V(t)>2\ve$ is very small.  Indeed,
\bel{bs}
\lim_{h\to 0+} ~{ \meas \bigl\{ t\in [\tau, \tau+h]\,;~~V(t)>2\ve\bigr\} \over h}~=~0.\eeq
%To understand the heart of the matter,
%let us compare the solutions to the two quasilinear hyperbolic equations
%$$\bega{cl} u_t+A(u) u_x&=~0,\\[1mm]
%w_t + A\bigl(u(\tau,y_k)\bigr)w_x&=~0,\enda\qquad\qquad w(\tau,x)= u(\tau,x)\quad
%\hbox{for}\quad x\in \,]x_{k-1}, x_k[\,.$$
At time $\tau+h$ the difference between $u$ and the solution 
$W_k$ to the linearized equation (\ref{linp}) can be estimated as
\bel{liner}
\bega{l}\ds E_k(h)~\doteq~\int_{x_{k-1}+h}^{x_k-h} \bigl| u(\tau+h,x) - W_k(\tau+h,x)\bigr|\, dx\\[4mm]
\qquad\ds
=~\O(1) \cdot \int_\tau^{\tau+h} \tv\bigl\{ u(t,\cdot);\, J_k(t)\bigr\} \cdot
\Big\| Df\bigl(u(t,\cdot)\bigr) -Df\bigl(u(\tau, y_k)\bigr)\Big\|_{\L^\infty}\, dt
\\[4mm]
\qquad\ds
=~\O(1) \cdot \int_\tau^{\tau+h} V(t)  \cdot
\bigl\| u(t,\cdot) -u(\tau, y_k)\bigr\|_{\L^\infty}\, dt.
%~=~\O(1)\cdot h\, \ve^2. 
\enda\eeq
If $V(t)\leq 2\ve$ for all $t\in [\tau, \tau+h]$, we would be in 
the Tame Variation case.  Both  factors in the integrand on the right hand side of 
(\ref{liner}) have size $\O(1)\cdot \ve$.
Hence $E_k=\O(1)\cdot h \ve^2$, as 
proved in  \cite{Bbook, BGlimm, BLF}.
In the general case,  there is an additional error, 
measured by  how much the solution $u$ can change
during the intervals of time where $V(t)>2\ve$.
If $\tau$ is a Lebesgue point for $V$, since the map $t\mapsto u(t,\cdot)$ is 
Lipschitz continuous,
by (\ref{bs}) we obtain the slightly weaker bound
$$E_k(h)~=~\O(1)\cdot
h\ve^2 + o(h),
$$
where the Landau symbol $o(h)$ denotes a higher order infinitesimal as $h\to 0$.
We now state a more precise result in this direction. As before, we assume that 
all characteristic speeds are contained in the interval $[-1,1]$.
\begin{lemma}\label{l:63}
Let $t\mapsto u(t,\cdot)$ be a weak solution to the strictly hyperbolic system (\ref{hcl}),
Lipschitz continuous with values in $\L^1(\R;\,\R^n)$. Assume that
\begi
\item[(i)]
$\tv\bigl\{ u(t,\cdot)\bigr\}\leq M$ for all $t\geq 0$.
\item[(ii)]
$\tv\bigl\{ u(\tau,\cdot)\,;~]a,b[\,\bigr\}\leq \ve$.
 \item[(iii)]
Setting
\bel{Vdef}V(t)~=~\tv\Big\{ u(t,\cdot)\,;~\bigl]a+(t-\tau)\,,~b-(t-\tau)\bigr[\Big\},\eeq
the limit (\ref{bs}) holds.
\endi
Call $W=W(t,x)$ the solution to the linear Cauchy problem with constant coefficients
\bel{linw} w_t +\Tilde A w_x~=~0,\qquad\qquad w(\tau,x) = u(\tau,x),\qquad x\in\R
%\,\quad\hbox{for}\quad x\in [a,b]
,\eeq 
with $\Tilde A= Df\bigl(u(\tau, \xi)\bigr)$ for some $\xi\in\,]a,b[\,$.
Then there holds
\bel{err6}\limsup_{h\to 0+} {1\over h} \int_{a+h}^{b-h}
\Big| u(\tau+h,x) - W(\tau+h,x)\Big|dx~=~\O(1)\cdot\ve^2.\eeq
\end{lemma}

{\bf Proof.}  By Theorem~4.3.1 in \cite{Dbook}, the bound (i) on the total variation
implies that the map $t\mapsto u(t,\cdot)$ is Lipschitz continuous with values in 
$\L^1(\R;\, \R^n)$.

Let $\lambda_i(u)$, $l_i(u)$, $r_i(u)$, $i=1,\ldots,n$,   be
respectively the eigenvalues and the  left and right eigenvectors of
the matrix $ A(u) = Df(u)$. For notational convenience, call 
$\tilde u = u(\tau,\xi)$, and let 
$$\tilde \lambda_i = \lambda_i(\tilde u),\qquad  \tilde l_i = l_i(\tilde u),\qquad 
 \tilde r_i= r_i(\tilde u),$$
be the corresponding eigenvalues and left and right eigenvectors of the matrix 
$\Tilde A = Df(\tilde u)$.

\begin{figure}[htbp]
\centering
 \includegraphics[scale=0.4]{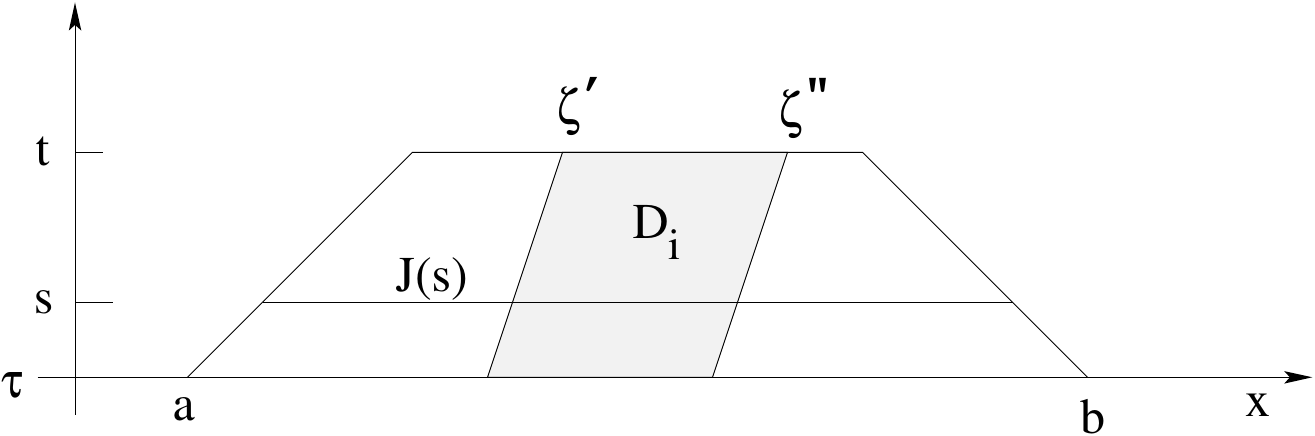}
    \caption{\small  The domain $D_i$ considered at (\ref{Di}).}
\label{f:claw85}
\end{figure}

Since $W$ solves the linear problem (\ref{linw}), one has
$$\tilde l_i\cdot W(t,x)~=~\tilde l_i\cdot W\big(\tau,
x-(t-\tau)\tilde\lambda_i\big)~=~\tilde l_i\cdot u\big(\tau,
x-(t-\tau)\tilde\lambda_i\big).$$
Following the proof of Theorem~9.4 in \cite{Bbook}, fix any two points
\bel{Jtdef}\zeta',\zeta''~\in~ J(t)~\doteq~\bigl] a+(t-\tau), \, b-(t-\tau)\bigr[\,.\eeq
Assuming $\zeta'<\zeta''$, consider the quantity
\bel{53}\bega{rl} E_i(\zeta',\zeta'')&\ds
\doteq ~ \tilde  l_i\cdot\int_{\zeta'}^{\zeta''}\bigl(u(t,x)
-W(t,x)\bigr)~dx\\[3mm]
&=~ \ds \tilde  l_i\cdot\int_{\zeta'}^{\zeta''}\Big(u(t,x)
-u\bigl(\tau,~x-(t-\tau)\tilde\lambda_i\bigr)\Big)~dx
.\enda \eeq
We apply the divergence theorem to the vector
$\big(u,\,f(u)\big)$
on the domain
\bel{Di}D_i~\doteq~\Big\{ (s,x);~~s\in[\tau,~t],~~
\zeta'-(t-s)\tilde\lambda_i\leq x\leq \zeta''-(t-s)\tilde
\lambda_i\Big\},\eeq
shown in Fig.~\ref{f:claw85}.

Since $u$ satisfies the conservation equation (\ref{hcl}),
the difference between the integral of $u$ at the top and at the
bottom of the domain $D_i$ is
measured by the inflow from the left side minus
the outflow from the right side of $D_i$.
By (\ref{53}) it thus follows
\bel{54}\bega{rl} E_i(\zeta',\zeta'')=& ~\ds
\int_\tau^{t} \tilde l_i\cdot\Big(\big(f(u)-\tilde\lambda_i 
u\big)(s,~\zeta'-(t-s)\tilde \lambda_i)\Big) \,ds\cr
&\ds \qquad-\int_\tau^{t} \tilde l_i\cdot\Big(\big(f(u)-\lambda_i 
u\big)(s,~\zeta''-(t-s)\tilde \lambda_i)\Big)\,ds\cr
=&\ds \int_\tau^{t}  l_i(\tilde u)\cdot\Big(\big(f(u'(s))-
\lambda_i\left(\tilde u\right) 
u'(s)\big)-  \big(f(u''(s))-\lambda_i(\tilde u) 
u''(s)\big)
\Big)\,ds\cr
=&\ds \int_\tau^{t}   F\left(\tilde u,u'(s),u''(s)\right)\,ds,
\enda \eeq
where we set
$$u'(s)\doteq u\big(s,~\zeta'-(t-s)\tilde\lambda_i\big),\qquad
u''(s)\doteq u\big(s,~\zeta''-(t-s)\tilde\lambda_i\big),$$
\begin{displaymath}
   F\left( u,u_{1},u_{2}\right)
  \,\doteq\,l_i( u)\cdot\Big(\big(f(u_{1})-\lambda_i(u) 
u_{1}\big)-  \big(f(u_{2})- \lambda_i (u)
u_{2}\big)
\Big).
\end{displaymath}
Observing that
\begin{itemize}
\item
  $ F\left(u,u_{2},u_{2}\right)\,=\,0$,
\item
  $D_{u_{1}} F\left(u,u_{1},u_{2}\right)\,=\,l_{i}(u)\cdot
  \left(Df(u_{1})-\lambda_{i}(u)I\right)$,
\item
  $D_{u_{1}} F\left(u,u,u_{2}\right)\,=\,l_{i}(u)\cdot
  \left(Df(u)-\lambda_{i}(u)I\right)\,=\,0$,
\end{itemize}
one obtains 
$$
  \bega{rl}    F\left(u,u_{1},u_{2}\right)&=~
     F\left(u,u_{1},u_{2}\right)-
     F\left(u,u_{2},u_{2}\right)\\[4mm]
    &=~\ds
    \int_{0}^{1}D_{u_{1}} F\left(u,u_{2}+\sigma \left(u_{1}-u_{2}
      \right),u_{2}\right)\; d \sigma \cdot \left(u_{1}-u_{2}\right)\\[4mm]
    &=~\ds
    \int_{0}^{1}\left[D_{u_{1}} F\left(u,u_{2}+\sigma \left(u_{1}-u_{2}
      \right),u_{2}\right)-D_{u_{1}} F\left(u,u,u_{2}\right)\right]\; d
  \sigma \cdot \left(u_{1}-u_{2}\right)\\[4mm]
  &=~\ds \O(1)\cdot \bigl(\left|u_{1}-u\right|+\left|u_2-u\right|\bigr)\cdot\left|u_{1}-u_{2}\right|.
  \enda $$
In turn, this yields
\bel{Eies}
    E_i(\zeta',\zeta'')~=~\O (1)\cdot \int_\tau^{t} \left|u'(s)-u''(s)\right|\cdot
    \Big(\left|u'(s)-\tilde u\right|+ \left|u''(s)-\tilde
        u\right|\Big)\, ds.
\eeq
Recalling the definitions at  (\ref{Vdef}) and (\ref{Jtdef}), 
 for any $x\in J(s)$ we now compute
\begin{displaymath}
    \left|u'(s)-\tilde u\right|\, \le\, V(s)+
    \left|u\big(s,x\big)- u\left(\tau,x\right) \right| + \left| u\left(\tau,x\right)-\tilde u\right|
  \,  \le\, V(s)+
    \left|u\big(s,x\big)- u\left(\tau,x\right) \right| + V(\tau).
\end{displaymath}
Integrating w.r.t.~$x$ over the interval
$J(s)$, dividing by its length and using the Lipschitz continuity of the map $t\mapsto u(t,\cdot)$,
 we obtain
\bel{gdef}\bega{rl}
    \left|u'(s)-\tilde u\right| &\ds \le~ V(s)+V(\tau) +  \frac{1}{\meas \bigl(J(s)\bigr)}\int_{J(s)}\bigl|u(s,x)-u(\tau,x)\bigr|\, dx\\[4mm]
    &\ds =~V(s)+\ve +
    \O(1)\cdot (s-\tau)~\dot=~g(s).
\enda \eeq
An entirely similar estimate holds for $\bigl|u''(s)-\tilde u\bigr|$. 
Therefore
$$
  \bega{rl}
    E_i(\zeta',\zeta'')&\ds =~\O (1)\cdot \int_\tau^{t} \left|u'(s)-u''(s)\right|\cdot g(s)\,
    ds\\[4mm]
    &=~\ds \O(1)\cdot  \int_\tau^t \tv\Big\{ u(s)\,;~\bigl]\zeta'-(t-s)\tilde\lambda_{i},~
    \zeta''-(t-s)\tilde\lambda_{i}\bigr]\Big\}\cdot
    g(s)\, ds\\[4mm]
   &=~\O(1)\cdot\mu_{i}\bigl(]\zeta',\zeta'']\bigr).
\enda
$$
Here $\mu_i$ is the Borel measure defined by
$$
  \mu_{i}\bigl(]c,d[\bigr)~=~  \int_\tau^t \tv\Big\{ u(s)\,;~\bigl]c-(t-s)\tilde\lambda_{i},~
   d-(t-s)\tilde\lambda_{i}\bigr[\Big\}\cdot
    g(s)\, ds, $$   for any open interval $]c,d[ \,\subset J(t)$.  
    
As proved in  Lemma~9.3 of \cite{Bbook}, one has
\begin{displaymath}
  \begin{split}
    \int_{J(t)}\bigl|u(t,x)-W(t,x)\bigr|\; dx&=~ \O(1)\cdot
    \sum_{i=1}^{n} \int_{J(t)}\left|\tilde l_{i}\cdot
      \bigl(u(t,x)-W(t,x)\bigr)\right|\; dx\\
    &=~
    \O(1)\cdot \sum_{i=1}^{n}\mu_{i}\left(J(t)\right)~=~\O(1)\cdot
    \int_{\tau}^{t}V(s)\cdot g(s)\;
    ds.
  \end{split}
\end{displaymath}
In turn, this implies
\bel{uuf}
\frac{1}{t-\tau}  \int_{J(t)}\bigl|u(t,x)-W(t,x)\bigr|\,
dx~=~ \O(1)\cdot\left(\frac{\left\|g\right\|_{\infty}}{t-\tau}\int_{\tau}^{t}\bigl|V(s)-V(\tau)\bigr|\, ds
+\frac{V(\tau)}{t-\tau}\int_{\tau}^{t}g(s)\, ds\right).
\eeq
We now observe that, for all $s>\tau$ sufficiently close to $\tau$,
the function $g$ introduced at (\ref{gdef}) satisfies
\bel{gbound}
g(s)~\leq~V(s) +\ve + \O(1) \cdot \left(t-\tau\right).\eeq
Thanks to (\ref{bs}), taking the limit of (\ref{uuf}) as $t  \to\tau+$ we thus obtain
\bel{limuuf}
  \limsup_{t\to\tau+}
  \frac{1}{t-\tau}  \int_{J(t)}\bigl|u(t,x)-W(t,x)\bigr|\;
dx~=~ \O(1)\cdot V(\tau)\left(V(\tau)+\varepsilon\right)~=~ \O(1)\cdot\varepsilon^{2},
\eeq
proving (\ref{err6}).
\endproof
\v

Going back to instantaneous error estimate (\ref{er1}), by Lemma~\ref{l:62}
it follows
\bel{Ak} \lim_{h\to 0+}~A_k(h)~=~0\qquad\qquad
\forall k=1,\ldots, N-1.\eeq
Moreover, using Lemma~\ref{l:63}, for some constant $C$
and every $k=1,\ldots, N$ we obtain
\bel{Bk} \limsup_{h\to 0+}~B_k(h)~\leq ~C\, \ve^2.\eeq
Observing that the number of intervals in the partition is $N = \O(1)\cdot \ve^{-1}$, 
we conclude
\bel{lim3}\limsup_{h\to 0+}\left(
 \sum_{k=1}^{N-1}  A_k(h)+  \sum_{k=1}^N B_k(h)\right)~\leq~0 + N \cdot C\ve^2 
 ~=~\O(1)\cdot \ve.\eeq
 Since $\ve>0$ can be chosen arbitrarily small, this implies (\ref{linf}).
\v
We can now state the main uniqueness theorem in \cite{BDL}.
\begin{theorem}\label{t:uniq}
Let (\ref{hcl}) be a strictly hyperbolic $n\times n$ system of conservation laws 
and consider the semigroup  of vanishing viscosity solutions $S:\D\times \R_+\mapsto \D$, constructed in Theorem~\ref{t:52}. Then, any  weak solution
$t\mapsto u(t,\cdot)\in \D$, which takes values in the domain of the semigroup and 
satisfies the Liu admissibility conditions (\ref{liuadm})  at every point of approximate jump,
coincides with a semigroup trajectory: 
\bel{uniq}
u(t,\cdot) ~=~ S_t u(0)\qquad\qquad \forall t\geq 0.\eeq
\end{theorem}

Notice that, by Theorem~\ref{t:52}, every limit of the vanishing viscosity approximations 
(\ref{evisc})  is a weak 
solution to (\ref{CP}) and satisfies the Liu admissibility conditions.
The above result provides a converse: every weak solution to (\ref{CP})
which is Liu-admissible
(and has suitably small total variation, so it lies within  the domain $\D$)
actually coincides with a semigroup trajectory. Therefore it is obtained as the 
unique limit of the viscous
approximations (\ref{evisc}).

\begin{figure}[htbp]
\centering
 \includegraphics[scale=0.5]{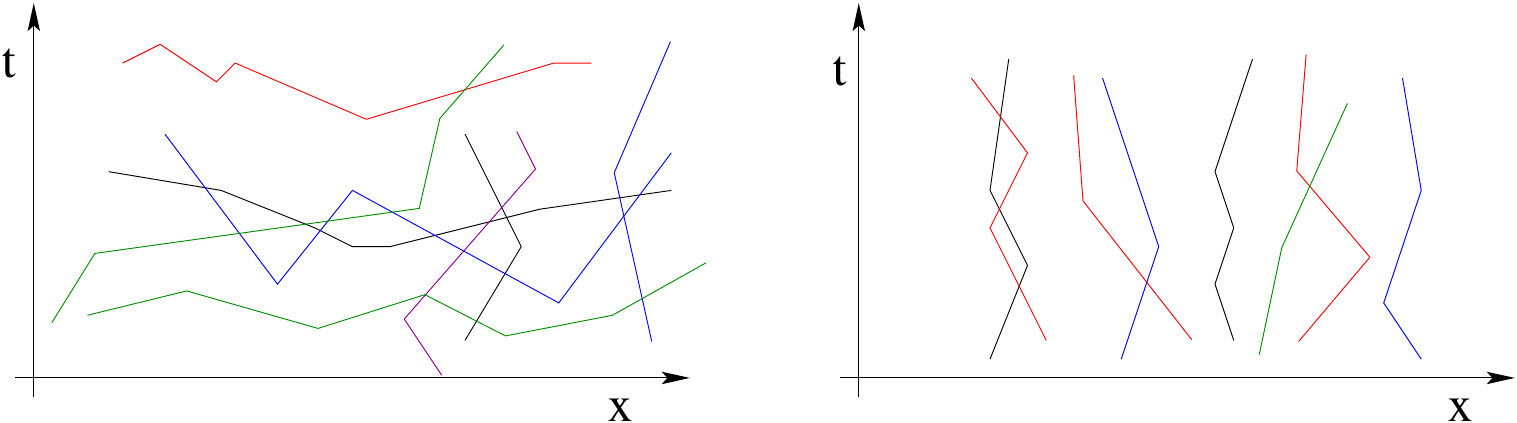}
    \caption{\small Left:  for a general $BV$ function 
    $u:\R_+\times\R\mapsto\R^n$ the set of approximate jumps can be covered by countably many Lipschitz curves.
    Right: if $u$ is a solution to (\ref{hcl}), all of its jumps travel with bounded speed.
    Hence the set of jumps can be covered by graphs of uniformly Lipschitz functions.}
\label{f:claw69}
\end{figure}

{\bf Sketch of the proof.} {\bf 1.} 
The assumption that $u=u(t,x)$ is a weak solution and its
total variation
$\tv\bigl\{ u(t,\cdot)\bigr\}$ remains uniformly bounded implies that $u:\R_+\times\R\mapsto\R^n$
is a BV function of the two variables $t,x$.   By a general structure theorem,
the set of its approximate jumps is countably rectifiable \cite{AFP}, i.e., it 
can be covered by 
countably many Lipschitz curves
(see Fig.~\ref{f:claw69}, left).   However, since $u$ is a solution to a hyperbolic system,
at each point of approximate jump the Rankine-Hugoniot equations hold.   In particular,
the speed of these jumps must be uniformly bounded.  
As proved in \cite{BDL}, the set of approximate jumps is contained in the graphs of 
countably many Lipschitz functions (see Fig.~\ref{f:claw69}, left). 
\bel{phil}x~=~\phi_\ell(t),\qquad\qquad \ell\geq 1.\eeq
To simplify our notation, w.l.o.g.~we shall 
assume that all characteristic speeds $\lambda_i(u)$ 
are contained in the interval
$[-1,1]$ and all functions $\phi_\ell$ have Lipschitz constant $1$.
\v
{\bf 2.} In addition to the functions $\phi_\ell$ we consider the countably many functions
\bel{phipm}\phi^{\xi+}(t)~=~\xi+t,\qquad\qquad \phi^{\zeta-}(t)~=~\xi-t,\eeq
where $\xi\in{\mathbb Q}$ is rational.  We relabel the set of all these
functions as
\bel{cmf}
\{ \phi_\ell\,;~\ell\geq 1\}\cup \{\phi^{\xi+}\,;~\xi\in{\mathbb Q}\} \cup \{\phi^{\zeta-}\,;~\zeta\in{\mathbb Q}\} ~=~\{\psi_j\,;~j\ge 1\}, \eeq
and consider the countably many functions
\bel{Wij}W_{ij}(t)~\doteq~\left\{\bega{cl} \tv\Big\{ u(t,\cdot)\,;~~\bigl] \psi_i(t),\,\psi_j
\bigr[\,\Big\}&\qquad \hbox{if}~~ \psi_i(t)<\psi_j(t),\\[1mm]
0 &\qquad \hbox{otherwise.}\enda\right.\eeq
We observe that each function $W_{ij}$ is measurable. 
Therefore there exists a null set $\N$ such that every 
$\tau\in \R_+\setminus \N$ is a Lebesgue point  for
all the countably many functions $W_{ij}$.
\v

\begin{figure}[htbp]
\centering
 \includegraphics[scale=0.43]{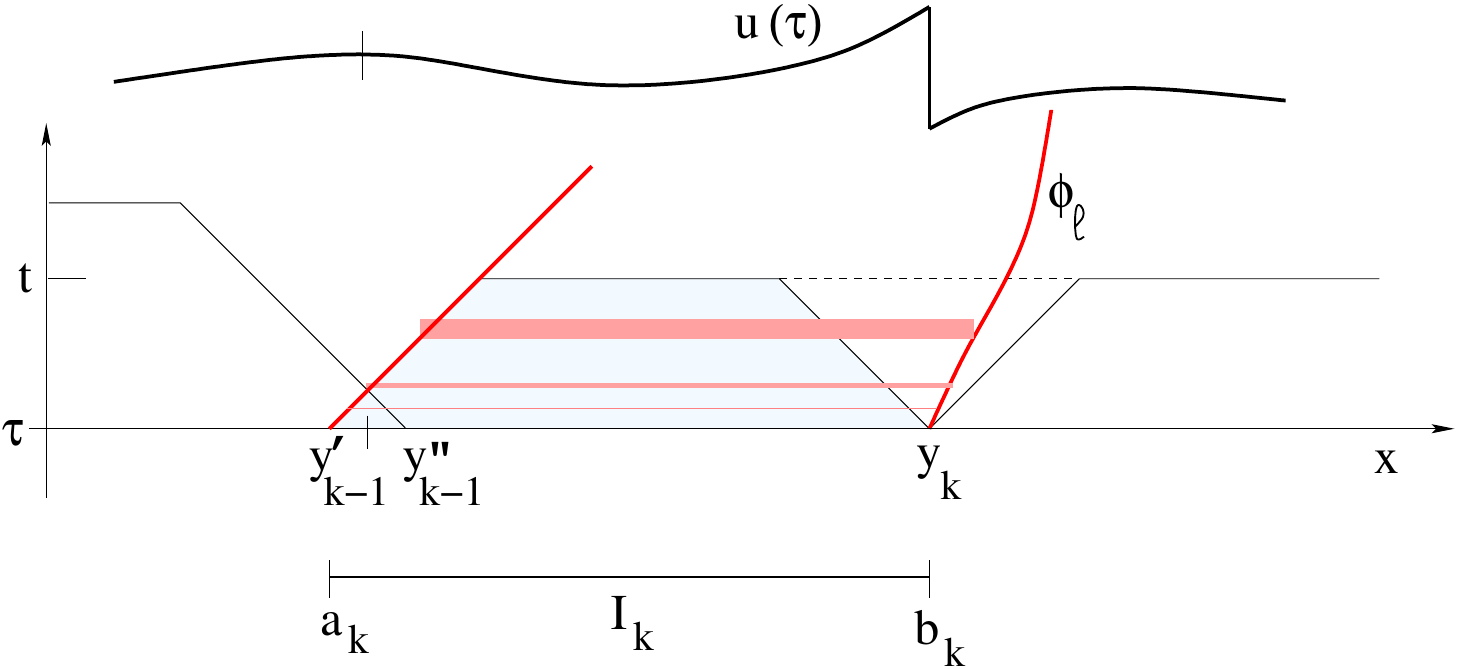}
    \caption{\small Covering the real line with points $y_k$ where $u(\tau,\cdot)$
    has a jump, and open intervals 
    $I_k= \,]a_k, b_k[\,$ where the total variation is $<2\ve$.}
\label{f:claw86}
\end{figure}

{\bf 3.} As shown in Fig.~\ref{f:claw86}, we now insert points $y_1<y_2<\cdots< y_N$
so that the total variation of $u(\tau,\cdot)$ on each open interval $\, ]y_{k-1}, y_k[\,$
is $<\ve$.   
\begi
\item
If $u(\tau,\cdot)$ has a jump at $y_k$, then by construction $(\tau, y_k)$ lies on the graph of one of the functions $\phi_\ell$.
\item If $u(\tau,\cdot)$  is continuous at $y_{k}$, then we can find two nearby
points $y'_k<y_k< y''_k$, lying on one of the rational lines
in (\ref{phipm}), namely
$$y_k'~=~\phi^{\xi+}(\tau),\qquad\qquad  y_k'' = \phi^{\zeta-}(\tau), \qquad
\hbox{for some}~~\xi,\zeta\in{\mathbb Q},$$
and furthermore
$$\tv\bigl\{ u(\tau,\cdot)\,;~[y_k', y_k'']\bigr\}~<~\ve.$$
\endi

In the end, we can cover the real line with finitely many points $y_k$ and 
open intervals $I_k=\, ]a_k, b_k[\,$ with the following properties.
\begi 
\item[(i)] At each point $y_k$, the function 
$u(\tau,\cdot)$ has an approximate jump, satisfying the Liu
admissibility condition.   
By Lemma~\ref{l:62} this implies
\bel{er22}\lim_{h\to 0+}  {1\over h} \int_{y_k-h}^{y_k+h}
 \bigg|u(t,x) - U(t-\tau,x-y_k)\bigg| \, dx ~=~0.\eeq

\item[(ii)]  On each open interval $I_k$ the total variation satisfies
$\tv\bigl\{ u(\tau,\cdot)\,;~I_k\bigr\}<2\ve$.   Moreover, both endpoints
of $I_k$  lie on the graph of one of the functions $\psi_j$
at (\ref{cmf}), say
$$a_k=\psi_i(\tau),\qquad \quad b_k = \psi_j(\tau),$$
for some $i,j\geq 1$.
Since all functions $\phi_\ell$ have Lipschitz constant $\leq 1$, for $t\geq \tau$ we have
$$\tv\Big\{u(t,\cdot)\,;~~\bigl]a_k+(t-\tau),\, b_k-(t-\tau)\bigr[\,\Big\}~\leq~\tv\Big\{u(t,\cdot)\,;~~\bigl]\psi_i(t),\, \psi_j(t)\bigr[\,\Big\}.$$
Therefore, calling $W_k$ the solution to the linearized equation
$$w_t+ \Tilde A_k w_x~=~0,\qquad\qquad  w(\tau,x)=u(\tau,x),\qquad \Tilde A_k\doteq Df\left(u\Big(\tau, {a_k+b_k\over 2}\Big)\right),$$
by Lemma~\ref{l:63} it follows
\bel{er66}\limsup_{h\to 0+} {1\over h} \int_{a_k+h}^{b_k-h}
\Big| u(\tau+h,x) - W_k(\tau+h,x)\Big|dx~=~\O(1)\cdot\ve^2.\eeq
\endi

Since 
 the semigroup trajectory $v(t,\cdot) =S_{t-\tau} u(\tau)$ satisfies the same estimates 
 (\ref{er22})-(\ref{er66}) as $u(t,\cdot)$,
we conclude that, 
at any Lebesgue time  $\tau\in \R_+\setminus \N$, 
\bel{lisup}
\limsup_{h\to 0+}~{\big\| u(\tau+h)-S_hu(\tau)\big\|_{\L^1}
\over h}~\leq~\sum_{k=1}^N  C_1\, \ve^2~\leq~C_2\,\ve,\eeq
for some constants $C_1, C_2$.
Indeed, the number of intervals $I_k$ in the partition is $N=\O(1)\cdot \ve^{-1}$.
Since $\ve>0$ can be chosen arbitrarily small, this achieves the proof.
\endproof

\section{Error estimates}
\label{sec:7}
\setcounter{equation}{0}
Having constructed a Lipschitz semigroup of admissible solutions to (\ref{hcl}),
it is of interest to estimate the $\L^1$ distance between an approximate solution
constructed by one of the algorithms described in Section~\ref{sec:3} 
and the exact solution.

   Results in this direction were proved in \cite{BGlimm}
for the front tracking method, in \cite{AM11, BM, BiMo} for the Glimm scheme, and in 
\cite{BHWY, BY} for vanishing viscosity approximations.

More precisely, consider the $n\times n$ 
hyperbolic system (\ref{hcl}), assuming that  all characteristic fields
are genuinely nonlinear.  
The estimate in \cite{BY} shows that the distance between the solution $u^\ve$
of the viscous approximation (\ref{diffa}) and the exact solution $u(t)= S_t \bar u$ 
 to the Cauchy problem
(\ref{CP3}) can be estimated as 
\bel{vap}\bigl\| u^\ve(t,\cdot)- u (t,\cdot)\bigr\|_{\L^1}~=~\O(1)\cdot (1+t) \, \tv\{\bar u\}\,
\sqrt\ve |\ln \ve|.\eeq
Next, consider  an approximate solution $u^{Glimm}$ constructed by the Glimm scheme, with a grid of step size
$\Delta t = \Delta x =\ve$.  Choosing sampling points as in (\ref{rap}),
the analysis in \cite{BM} has established a similar convergence rate:
\bel{Glimmes}\lim_{\ve\to 0} {\bigl\| u^{Glimm} (t,\cdot) - u(t,\cdot)\bigr\|_{\L^1}\over 
\sqrt\ve |\ln \ve|}~=~0\qquad\qquad \forall t> 0.\eeq

For other approximation methods, such as periodic mollifications, the
backward Euler scheme, or fully discrete numerical schemes,
no a priori BV bounds are currently available.  
In particular, it is known that the Godunov scheme can amplify the total 
variation by an arbitrarily large factor~\cite{BBJ}.
For this reason, it seems more promising to look for {\it a posteriori} error bounds.
Namely, assume that an approximate solution to the Cauchy problem (\ref{CP3}) 
has been constructed, whose total variation
remains small for all times $t\in [0,T]$.    Using this additional information, we seek 
a bound on the error
\bel{error}
\bigl\| u^{approx}(t,\cdot) - u^{exact}(t,\cdot)\bigr\|_{\L^1}\,.\eeq
We outline here an approach which is in a sense ``universal", 
i.e., it does not make reference
to any particular approximation method.
%
%We seek a posteriori error estimates on approximate solutions 
%$u= u(t,x)$ to a Cauchy problem (\ref{CP3}).
%%\bel{CP2} u_t + f(u)_x ~=~0,\qquad\qquad u(0,x) = \bar u(x).\eeq
%As usual, we assume that this $n\times n$ system is strictly hyperbolic, and generates a Lipschitz
%semigroup of entropy-weak solutions $u(t) = S_t \bar u$ on a domain 
%$\D\subset\L^1(\R;\,\R^n)$ containing all functions with sufficiently small total variation 
%\cite{BiB, Bbook, BLY}.
%More precisely,  for some constant $L$, the semigroup $S:\D\times\R_+\mapsto \D$ 
%satisfies
%\bel{Slip}
%\bigl\| S_t\bar u - S_s\bar v\bigr\|_{\L^1}~\leq~L\Big( |t-s| + \|\bar u - \bar v\|_{\L^1}\Big)\qquad
%\forall s,t\geq 0,~~\bar u,\bar v\in\D.\eeq
%

Given $\ve>0$, consider an approximate solution $u=u(t,x)$ with the following properties.
\begin{definition}\label{d:as}
Let (\ref{hcl}) be an $n\times n$ strictly hyperbolic systems of conservation laws,
endowed with a strictly convex entropy $\eta$, with entropy flux $q$.
We say that $u=u(t,x)$ is an {\bf $\ve$-approximate solution} to the Cauchy problem
(\ref{CP}) if $\bigl\|u(0,\cdot)-\bar u\bigr\|_{\L^1}\leq\ve$ and moreover the following holds.
\begi\item[{ \bf (AL$_\ve$)}] {\bf Approximate Lipschitz continuity:}
$$\|u(\tau,\cdot)-u(\tau',\cdot)\|_{\L^1}~\leq~M\,|\tau-\tau'|+\ve
%\cdot\sup_{t\in [\tau, \tau']} \tv\bigl\{ u(t,\cdot)\bigr\}
\qquad
\qquad \forall \tau,\tau' \geq 0$$
\item[{\bf (P$_\ve$)}] {\bf Approximate conservation law and approximate entropy inequality:}

{\it For every strip $[\tau,\tau' ]\times \R$  and every test function $\vp\in \C^1_c(\R^2)$, one has}
\bel{wsol}\bega{l}\ds \left| \int u(\tau,x)\vp(\tau,x)\, dx-\int u(\tau',x)\vp(\tau,x)\, dx+
\int_\tau^{\tau'}\! \!\int \bigl\{ u\vp_t+f(u)\vp_x\bigr\}\,dx\,dt\right|\\[4mm]
\qquad\qquad\ds \leq~ \ve (\tau'-\tau+\ve)\,\|\vp\|_{W^{1,\infty}}\,
%\cdot \sup_{t\in [\tau, \tau']} \tv\bigl\{ u(t,\cdot)\bigr\}
.\enda\eeq

{\it
Moreover, assuming $\vp\geq 0$, one has the entropy inequality}
\bel{eiq}
\bega{l}\ds  \int \eta(u(\tau,x))\vp(\tau,x)\, dx-\int \eta(u(\tau',x))\vp(\tau',x)\, dx+
\int_\tau^{\tau'}\! \!\int \bigl\{\eta(u)\vp_t+q(u)\vp_x\bigr\}\,  dxdt
\\[4mm]\qquad\qquad\ds \geq~- \ve(\tau'-\tau+\ve)\, \|\vp\|_{W^{1,\infty}}
%\cdot \sup_{t\in [\tau, \tau']} \tv\bigl\{ u(t,\cdot)\bigr\}
.\enda\eeq
\endi
\end{definition}

\begin{figure}[htbp]
\centering
  \includegraphics[scale=0.4]{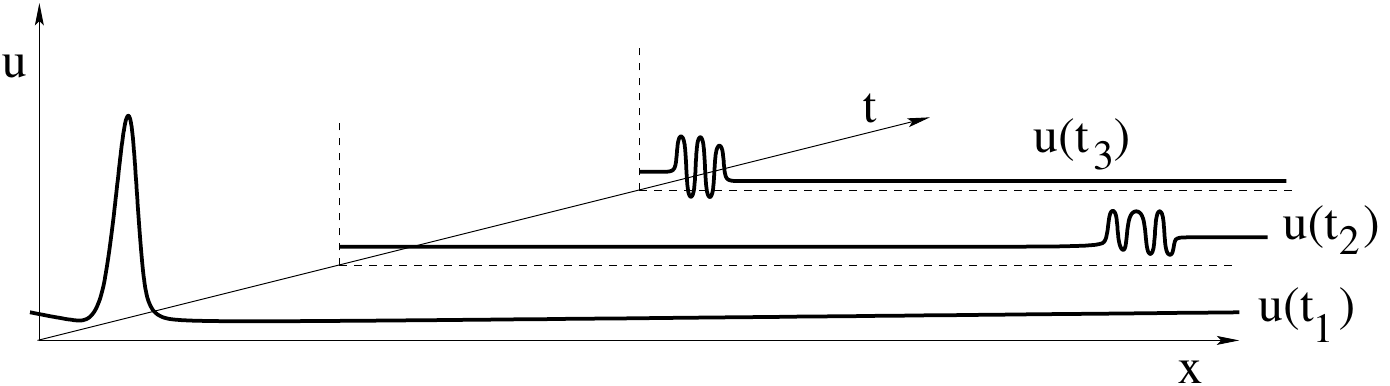}
    \caption{\small An approximate solution $u$ 
    where the total variation remains small at all times. However,
   oscillations appear and disappear at different regions on the $x$-$t$ plane.  }
\label{f:hyp221}
\end{figure}

In the above setting, the paper \cite{BCS} has established {\bf a posteriori}
error estimates,
assuming that the total variation of the $\ve$-approximate 
solution remains small, so that $u(t,\cdot)$ remains within the domain of the semigroup.  However, the estimates in \cite{BCS} also required a
``post processing algorithm", tracing the location of the large shocks in the approximate solution.
This is related to the assumptions of  ``tame variation", ``tame oscillation" or ``bounded variation along space-like curves"
which were used respectively in \cite{BLF}, \cite{BG2} and in \cite{BL} to prove uniqueness of solutions.  In essence, these additional 
assumptions rule out configurations such as the one shown in Fig.~\ref{f:hyp221}.

The recent paper \cite{BGu} has shown that, for a system endowed with a strictly convex entropy, these additional regularity conditions are not needed to achieve uniqueness:

\begin{theorem}\label{t:1} Let (\ref{hcl}) be a strictly hyperbolic $n\times n$ system, where
each characteristic field is either genuinely nonlinear or linearly degenerate, and which
admits a strictly convex entropy $\eta(\cdot)$. 
Then every entropy admissible weak solution $u:[0,T]\mapsto\D$, coincides with a semigroup trajectory.  \end{theorem}

We observe that, in the setting of the above theorem, the dissipation of a single entropy
suffices to single out the Liu-admissible shocks.

As proved in \cite{BGu}, the compactness of the family of approximate solutions, 
together with the uniqueness of the limit, yields a uniform convergence rate:
\begin{corollary}
In the above setting, given $T,R>0$, there exists a function $\ve\mapsto \varrho(\ve)$ 
with the following properties.
\begi
\item[(i)]
$\varrho$ is continuous, nondecreasing, with $\varrho(0)=0$. 
\item[(ii)] Let $t\mapsto u(t,\cdot)\in \D$ be an $\ve$-approximate solution to 
(\ref{CP3}), with $u(t,\cdot)$ supported inside the interval $[-R,R]$ for all $t\in [0,T]$.
Then one has
\bel{err4}
\bigl\|u(t)-S_t \bar u\bigr\|_{\L^1}~\leq~\varrho(\ve)\qquad\qquad \forall t\in [0,T].\eeq
\endi
\end{corollary}
This corollary shows that such a ``universal rate of convergence" must exist.
However, it does not offer clues on how  the function $\rho(\cdot)$ looks like.

\v
{\bf Open Problem \#3.} {\it Let (\ref{hcl}) be a strictly hyperbolic $n\times n$ system, where
each characteristic field is either genuinely nonlinear or linearly degenerate, and which
admits a strictly convex entropy $\eta(\cdot)$. Provide an asymptotic estimate
on the universal convergence rate $\varrho(\cdot)$ in (\ref{err4}), as $\ve\to 0$.}

Based on the earlier estimates (\ref{vap})-(\ref{Glimmes}), 
in the genuinely nonlinear case one may conjecture that 
$\varrho(\ve)\approx \ve^{1/2} |\ln\ve|$.    

The key feature of the bound (\ref{err4}) is that it holds for any $\ve$-approximate solution
satisfying {\bf (AL$_\ve$)-(P$_\ve$)}, regardless of the method used 
to construct the approximation.  All the algorithms considered in Section~\ref{sec:3}
generate $\ve$-approximate solutions, in the sense of Definition~\ref{d:as}.
See Section~6 in \cite{BCS} for details.

\v
One can speculate whether a similar universal convergence rate 
can be valid for general $n\times n$ systems, not necessarily endowed with a 
strictly convex entropy.  For these systems, semigroup trajectories are characterized by 
the Liu admissibility condition, as in Definition~\ref{d:liu}.
To reach our goal, we should replace the 
$\ve$-approximate entropy condition (\ref{eiq}) with some sort of $\ve$-approximate
Liu condition.  This leads to
\v
{\bf Open Problem \#4.} {\it Introduce a definition of ``$\ve$-approximate Liu admissible
solution", valid for general $n\times n$ hyperbolic systems, possibly not endowed with a strictly convex entropy.}

A bit more precisely, what is needed here is a suitable definition such that the following properties will be satisfied.
\begi
\item Approximate solutions with small total variation
constructed by the various methods  described in Section~\ref{sec:3} should 
all satisfy the $\ve$-approximate Liu condition, with $\ve\to 0$ as the 
step size in the approximation (or the viscosity coefficient) approaches zero.
\item Given a convergent sequence of approximations $u_n\to u$, if each $u_n$ 
is an $\ve_n$-approximate Liu admissible solution with $\ve_n\to 0$, 
then the limit solution $u$ should be
Liu-admissible in the original sense.
\endi

\section{Solutions with unbounded variation}
\label{sec:8}
\setcounter{equation}{0}
As remarked in Section~\ref{sec:4},  for solutions with large initial
data it is a hard open question to decide whether the total variation remains bounded for all times.  
It is thus natural to consider solutions in the larger space $\L^\infty(\R;\R^n)$, possibly with unbounded variation.   

For $2\times 2$ systems, existence of weak solutions with $\L^\infty$ data
was proved in a fundamental paper by DiPerna~\cite{DP83}, based on compensated compactness.  See also \cite{Dbook, Lu, Serre}
for a comprehensive account of this approach. 
Existence of $\L^\infty$ solutions  remains a largely open problem for general $n\times n$ systems.

Unfortunately, compensated compactness works as a ``black box".  It provides
an abstract result on the existence of solutions, but it does not yield information
about uniqueness, continuous dependence, or the qualitative structure 
of these solutions. Some of the few results on the regularity of
$\L^\infty$ solutions
can be found in \cite{CT, Golding}.

In this direction, it would be of interest to construct a continuous 
semigroup of admissible solutions, defined on a domain larger than $BV$.
\v
{\bf Open Problem \#5.} {\it Given an $n\times n$ hyperbolic system of 
conservation laws, extend the semigroup of vanishing viscosity solutions to a larger
domain $\Tilde \D\subseteq\L^\infty(\R;\,\R^n)$, also 
containing functions with unbounded variation.}
\v
A continuous semigroup of solutions defined on the entire space $\L^\infty(\R;\,\R^n)$ was constructed in \cite{BG2} 
for some Temple class systems, and  more recently in \cite{BGS} 
for $2\times 2$ systems in triangular form. But apart from a few special cases
the problem is wide open.

As suggested in \cite{ABBCN}, 
in general it may not be possible to construct a continuous  semigroup 
defined on the entire space $\L^\infty$.  Instead, one could consider some intermediate domain $\Tilde \D\subset\L^\infty$, 
borrowing ideas from the theory of intermediate spaces used in the analysis of
parabolic equations 
\cite{Henry, Lunardi}.
Of course, we do not expect that the extended semigroup will be Lipschitz continuous.  
 Its modulus of continuity will strongly depend on the regularity 
properties of 
functions $u\in \Tilde D$.   

In addition to compensated compactness, 
another approach is worth mentioning here.
In their classical memoir~\cite{GL}, Glimm and Lax consider the Cauchy problem
for a genuinely nonlinear  $2\times 2$ system. 
Assuming that the initial data $\bar u$ 
has sufficiently small $\L^\infty$ norm, they prove that a global weak solution exists, globally in time.   
Indeed, the total variation (which initially may well be infinite) becomes 
locally finite at every  time $t>0$.   
See also \cite{BCM} for a shorter proof, based on front tracking approximations.
The uniqueness and continuous dependence of these solutions still remains an open problem.
\v
In the opposite direction, it would also be of interest to find examples of Cauchy problems admitting multiple solutions. In~\cite{BS1} a $3\times 3$ 
strictly hyperbolic system has been constructed, together with bounded, measurable initial data, leading to an infinite number of solutions.  However, this example
does not have physical relevance because the system does not admit 
convex entropies.   We thus conclude with
\v
{\bf Open Problem \#6.} {\it Construct an example of an $n\times n$ strictly
hyperbolic system, endowed with a strictly convex entropy, together with
initial data $\bar u\in \L^\infty(\R;\,\R^n)$, such that the   Cauchy problem
admits two distinct entropy admissible solutions.}

\v
{\bf Acknowledgement.}  This research
 was partially supported by NSF with
grant  
  DMS-2306926, ``Regularity and approximation of solutions to conservation laws".

\vs

\end{document}